\documentclass[11pt]{amsart}

\usepackage{graphicx}

\newtheorem{theorem}{Theorem}
\newtheorem{proposition}[theorem]{Proposition}
\newtheorem{lemma}[theorem]{Lemma}
\newtheorem{corollary}[theorem]{Corollary}

\theoremstyle{definition}
\newtheorem{definition}[theorem]{Definition}

\theoremstyle{remark}
\newtheorem{remark}[theorem]{Remark}

\def\R{\mathbb{R}}
\def\Z{\mathbb{Z}}
\def\d{\partial}
\def\im{\operatorname{im}}
\def\A{\mathbb{A}}
\def\C{\mathcal{C}}
\def\eps{\varepsilon}
\def\jetO{j^1(0)}
\def\jetS{\mathcal{J}^1(S^1)}
\def\L{\mathcal{L}}

\numberwithin{theorem}{section}

\begin{document}

\title{Legendrian Solid-Torus Links}
\author[L. Ng]{Lenhard Ng}
\address{Department of Mathematics, Stanford University, Stanford, CA 94305}
\author[L. Traynor]{Lisa Traynor}
\address{Department of Mathematics, Bryn Mawr College, Bryn Mawr, PA 19010}
\thanks{
Ng is supported by a Five-Year Fellowship from the American
Institute of Mathematics. Traynor was supported in part by NSF grant
DMS 9971374.
}

\begin{abstract}
Differential graded algebra invariants are constructed for
Legendrian links in the $1$-jet space of the circle.  In parallel to
the theory for $\R^3$, Poincar\'e--Chekanov polynomials and
characteristic algebras can be associated to such links. The theory
is applied to distinguish various knots, as well as links that are
closures of Legendrian versions of rational tangles. For a large
number of two-component links, the Poincar\'e--Chekanov polynomials
agree with the polynomials defined through the theory of generating
functions. Examples are given of knots and links which differ by an
even number of horizontal flypes that have the same polynomials but
distinct characteristic algebras. Results obtainable from a
Legendrian satellite construction are compared to results obtainable
from the DGA and generating function techniques.
\end{abstract}

\maketitle

\section{Introduction}

In the late 1990's, there were some breakthrough ideas for
constructing new invariants of Legendrian links in $\R^3$. These
invariants came from adapting the techniques of pseudoholomorphic
curves, a powerful tool in the study of symplectic manifolds, to the
setting of contact manifolds, via contact homology \cite{bib:EGH}.
Combinatorial ways of calculating these ``holomorphic" invariants
were developed from both the Lagrangian and the front projections of
a Legendrian link \cite{bib:Ch,bib:ENS,bib:Ng1}.

In this paper, the focus will be to study
links in the 1-jet space of the circle, $\jetS$, a manifold diffeomorphic
to the solid torus $S^1 \times \R^2$:
$$\jetS= T^*(S^1) \times \R = \{ (x,y,z) \,:\, x \in S^1, \  y,z \in
\R \},$$
with contact structure given by $\xi = \ker( dz - y\,dx)$.
Viewing $S^1$ as a quotient of the unit interval, $S^1 =
[0,1]/(0\sim 1)$, we can visualize knots in $\jetS$ as quotients of
arcs in $I \times \R^2$ with appropriate boundary conditions.

Links in $\jetS$ were examined by  one of the authors  using the
technique of generating functions, \cite{bib:Tr}. In this paper
holomorphic techniques are developed to study these and other links.
In Section~\ref{sec:invariants}, the theory underlying the
holomorphic invariants is developed. Given the ``Lagrangian''
$xy$-projection of a Legendrian link in $\jetS$, it is possible to
combinatorially define the Chekanov--Eliashberg differential graded
algebra (DGA) over $\Z_2$, in a manner exactly following the
definition in $\R^3$ from \cite{bib:Ch}. However, as in $\R^3$, it
is often more convenient to work in the ``front'' $xz$-projection
rather than the Lagrangian $xy$-projection. To do this, we introduce
a suitable modification of the resolution technique for $\R^3$ from
\cite{bib:Ng1}, which produces a Lagrangian projection from a front
projection.
Using this resolution, and in parallel to the $\R^3$ theory, we can
associate to Legendrian links in $\jetS$ invariants such as
Poincar\'e--Chekanov polynomials, which measure the homology of a
linearized version of the Chekanov--Eliashberg DGA, and
characteristic algebras, which measure some nonlinear information
from the DGA. We also make use of the additional structure on the
DGA provided by Mishachev's homotopy splitting, \cite{bib:Mi}. For
multi-component links, this results in split Poincar\'e--Chekanov
polynomials and split characteristic algebras.

\begin{figure}[ht]
\centerline{
\includegraphics[height=1in]{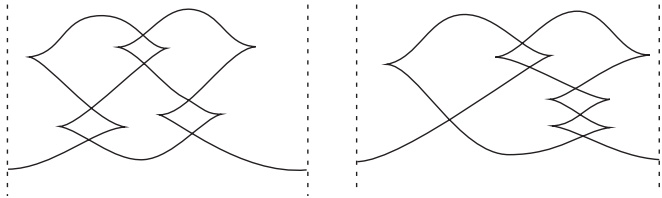}
}
\caption{
Closures of long knot versions of Chekanov's $5_2$ knots. These can
be distinguished by the Poincar\'e--Chekanov polynomials.
}
\label{fig:intro_long}
\end{figure}

Natural knot candidates for study with this technique are formed by
identifying the ends of long Legendrian knots in $\R^3$.
Figure~\ref{fig:intro_long} shows two distinct Legendrian knots that
are similar in spirit to Chekanov's $5_2$ examples.  These can be
distinguished by the Poincar\'e--Chekanov polynomials.  More
examples of this sort are explored in Section~\ref{sec:knots}.

Some natural link candidates for study are the families of links
studied in \cite{bib:Tr}.   Such links can be specified by a vector
of the form
\begin{align*}
(2h_n, v_{n-1}, 2h_{n-1}^{p_{n-1}}, \dots, 2h_2^{p_2}, v_1, 2h_1^{p_1}), \quad
&h_n, v_{n-1}, \dots, h_2, v_1 \geq
1,  \quad h_1 \geq 0, \\
\text{ and } p_i \in &\{ 0, \dots, 2h_i \} \text{ for } i \in \{ 1,
\dots, n-1\}.
\end{align*}
These links can be seen as quotients of Legendrian versions of  even
parity rational tangles defined by Conway in \cite{bib:Co}. The
superscripts denote a link obtained by ``horizontal flypes" of a
standard configuration. In Section~\ref{sec:scp}, the
Poincar\'e--Chekanov polynomials of these links are calculated.
There is a striking parallel between the split Poincar\'e--Chekanov
polynomials defined via the holomorphic technique and the signed
polynomials defined via the generating function technique:  after a
change of variables, the positive Poincar\'e--Chekanov polynomial
agrees with the negative generating function polynomial, see
Remark~\ref{rmk:holo_vs_gf}. This gives evidence for the following
principle, which seems to hold for other contact manifolds as well:
\vspace{11pt}

\noindent \textit{The following two Legendrian-isotopy invariants
contain the same information, when defined:
\begin{itemize}
\item
Morse-theoretic Poincar\'e polynomials obtained from generating
functions;
\item first-order Poincar\'e--Chekanov polynomials
obtained from holomorphic curves and contact homology.
\end{itemize}}

\vspace{11pt}

\noindent This principle has been verified for all solid-torus links
for which generating function polynomials have been defined. At
present, this does not include all solid-torus links, and so there
are instances in which the
linearized DGA is applicable but generating functions are not.
In addition, we
will present an example of two links for which the generating
function polynomials are defined and equal, but which can be
distinguished by the full DGA via characteristic algebras.

\begin{figure}[ht]
\centerline{
\includegraphics[height=2in]{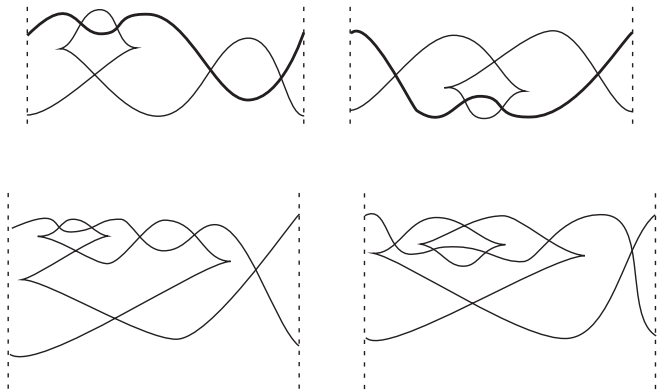}
}
\caption{
The first example gives the links $(2,1,2)$ and $(2,1,2^1)$; the
second example gives the knots $(2,1,2,1,1)$ and $(2,1,2^1,1,1)$.
These links and knots can be distinguished by their polynomials.
   }
\label{fig:link+knot_flypes}
\end{figure}

We now look at some concrete examples of links which can be
distinguished by one means or another. The first example in
Figure~\ref{fig:link+knot_flypes} gives topologically equivalent
links that do not have the same polynomials (generating function or
Poincar\'e--Chekanov) and are hence not Legendrian isotopic. The
second example gives knots of a similar flavor that can also be
distinguished by their polynomials. This knot example solves
Question 1.34 posed in \cite{bib:Tr}.

\begin{figure}[ht]
\centerline{
\includegraphics[width=4in]{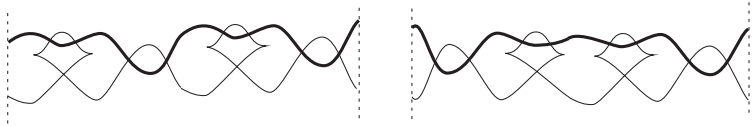}
}
\caption{
The Legendrian links $(2,1,2) \# (2,1,2)$ and  $(2,1,2^2) \#
(2,1,2)$, denoted in Section~\ref{sec:connectsums} as $\mathcal
L_{2,2}$ and $\mathcal L_{0,4}$. They have the same generating
function/Poincar\'e--Chekanov polynomials but different
characteristic algebras. }
\label{fig:intro_connect_sum}
\end{figure}

\begin{figure}[ht]
\centerline{
\includegraphics[width=4in]{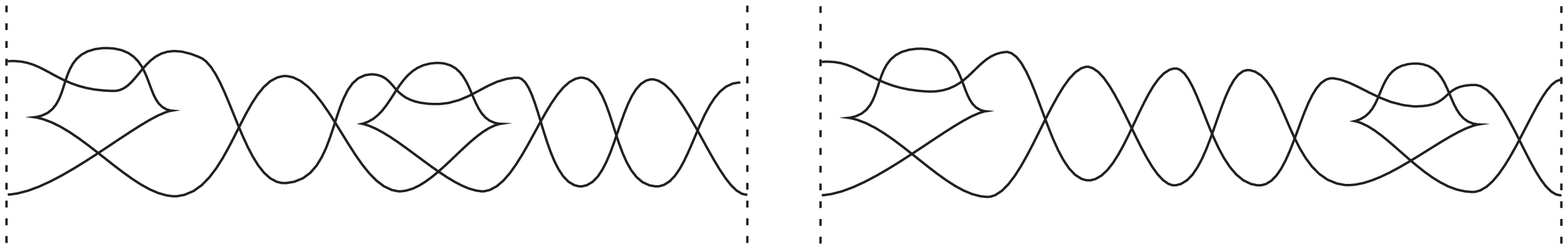}
}
\caption{ Two Legendrian knots quite similar to the examples in
Figure~\ref{fig:intro_connect_sum}, but with isomorphic algebras.
Are they isotopic?
}
\label{fig:L_2_3}
\end{figure}

Theorems~\ref{thm:StandardForm} and  \ref{thm:FlypeForm} show that
the polynomials may detect the parity
of the numbers of horizontal flypes.  This raises interesting
questions about links
that differ by an even number of flypes.  For example,
   Question 1.31 in \cite{bib:Tr} asks whether the Legendrian links denoted by
   $ (2,1,2) \# (2,1,2)$  and  $(2,1,2^2) \# (2,1,2)$ shown in
Figure~\ref{fig:intro_connect_sum} are isotopic.
In parallel to the generating function theory, the holomorphic
polynomials of these (nonrational) connect sums of the form $L_1 \#
L_2$ can be easily calculated from the polynomials of $L_1$ and
$L_2$, see Theorem~\ref{thm:connect_sum} and
Remark~\ref{rmk:gfconnectsum}.  Thus the polynomials of $(2,1,2) \#
(2,1,2)$  and $(2,1,2^2) \# (2,1,2)$  agree.  However, a calculation
of their characteristic algebras shows that these are in fact
distinct; see Proposition~\ref{prop:L_{k_1,k_2}_vs_L_{0,j_2}}. A
solution to this problem raises some new interesting questions. For
example, a slight modification of the  constructions in
Figure~\ref{fig:intro_connect_sum} produces the knots in
Figure~\ref{fig:L_2_3} that have tamely isomorphic DGAs.  Are these
knots isotopic? This and some related questions are posed in
Section~\ref{sec:connectsums}.

We include a short appendix which describes a technique, different
from the DGA or generating functions, which can also distinguish
Legendrian solid-torus links. The relevant construction, which we
call a Legendrian satellite, constructs a link in standard contact
$\R^3$ from a link in the solid torus. We show that Legendrian
satellites can be used to recover some of the results derived in
this paper using the DGA or in \cite{bib:Tr} using generating
functions, and can even be applied in cases where neither of the
other techniques works.

\subsection*{Acknowledgments}
We thank the American Institute of Mathematics for sponsoring the
contact geometry program during Fall 2000, where the ideas for this
paper originated.  We also thank the hospitality of the Institute
for Advanced Study in 2001--2002 when further ideas were developed.

\section{Holomorphic Invariants}
\label{sec:invariants}

\subsection{The DGA invariant}
\label{ssec:dgadef}

As in standard contact $\R^3$, apart from topological type, oriented
Legendrian links in $\jetS$ have two classical contact invariants,
rotation number $r$ and Thurston--Bennequin number $tb$. These are
defined in precisely the same way as in $\R^3$; for instance, in the
``front'' $xz$-projection, $2r$ is the number of cusps oriented
downwards minus the number oriented upwards, and $tb$ is the number
of crossings, counted with the usual signs, minus the number of
right cusps. In \cite{bib:Ch}, inspired by work of Eliashberg and
Hofer, Chekanov introduced a nonclassical invariant of Legendrian
knots in $\R^3$, and we will be interested in its analogue in
$\jetS$.

Given the ``Lagrangian" $xy$-projection of a Legendrian link in
$\jetS$, we may combinatorially define the Chekanov--Eliashberg
differential graded algebra $(\A,\d)$ in a manner exactly following
the definition in $\R^3$ from \cite{bib:Ch,bib:ENS}.  All proofs
translate, without changes, to solid-torus links.  We summarize here
the main results about the DGA; please refer to one of
\cite{bib:Ch,bib:ENS,bib:Ng1} for definitions of terms.

\begin{theorem} \label{thm:d2}
$\d^2=0$, and $\d$ lowers degree by $1$.
\end{theorem}

\begin{theorem} \label{thm:inv}
Two Legendrian-isotopic links in $\jetS$ have equivalent DGAs.
\end{theorem}

\noindent As will be explained in more detail below, all of the
algebraic machinery built around DGAs, including Poincar\'e
polynomials \cite{bib:Ch} and characteristic algebras
\cite{bib:Ng1}, can be applied to our setup. However, as in standard
contact $\R^3$, it is often more convenient to work in the front
projection rather than the Lagrangian projection.  The following
{\it resolution} technique, a slight modification of that used in
$\R^3$ in \cite{bib:Ng1}, tells us how to do the ``bookkeeping" when
working in the front projection. Using resolution, we will then
define the Chekanov--Eliashberg DGA for fronts of Legendrian links
in $\jetS$.

\subsection{Resolution}
\label{ssec:res}

\begin{figure}[ht]
\centerline{
\includegraphics[height=1in]{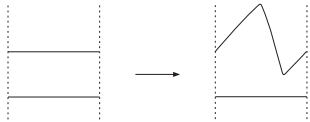}
}
\caption{
A modification of the $xz$-projection  as part of the resolution
procedure.  Since the upper strand is replaced by a segment of slope
$1$, a transition tangle segment must be inserted to guarantee that
the tangle closes to a link.
}
\label{fig:res0}
\end{figure}

\begin{figure}[ht]
\centerline{
\includegraphics[height=1in]{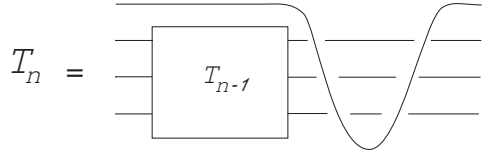}
}
\caption{
Inductive construction of the tangle $T_n$.
}
\label{fig:tangle}
\end{figure}

To convert an $xz$-projected Legendrian link to an $xy$-projected
link, we use resolution, as in \cite{bib:Ng1}. Figure~\ref{fig:res0}
illustrates for the link with an $xz$-projection consisting of two
line segments of slope zero
    how in the resolution
procedure the front diagram is modified so that initially the upper
strand has slope $1$.    As seen in this case, a complication not
present in the $\R^3$ case is that after the normal resolution
procedure, in order for the strands to close up smoothly under the
identification, a ``transition" tangle (with no crossings) needs to
    be adjoined to the $xz$-projection.
In the $xy$-projection, this translates into adjoining
a topologically trivial tangle $T_n$, with $n$ strands,
defined inductively as follows: $T_1$ is the trivial tangle
with $1$ strand, and $T_n$ is obtained from $T_{n-1}$ through the
construction shown in Figure~\ref{fig:tangle}.

\begin{figure}[ht]
\centerline{
\includegraphics[width=4in]{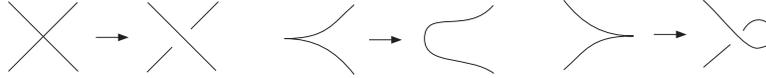}
}
\caption{
Resolving singularities in a front.
}
\label{fig:res}
\end{figure}

\begin{figure}[ht]
\centerline{
\includegraphics[height=2in]{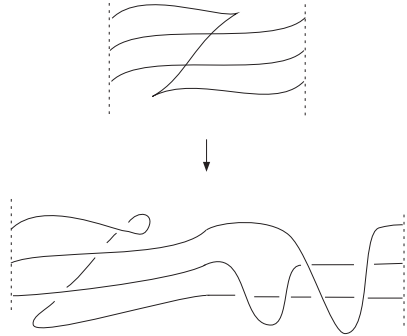}
}
\caption{
Example of a resolution.  The top diagram is a front projection; the
bottom diagram is a Lagrangian projection.
}
\label{fig:res-ex}
\end{figure}

\begin{definition} \label{def:resolution}
Let $Z$ be the front projection of a Legendrian link in  $S^1 \times
\R^2$. The {\it resolution} of $Z$ is the Lagrangian projection of a
Legendrian link obtained by resolving each singularity as in
Figure~\ref{fig:res}, and appending a tangle $T_n$ at the very right
of the projection, as given in Figure~\ref{fig:tangle}, for
appropriate $n$. See Figure~\ref{fig:res-ex} for an example.
\end{definition}

\begin{proposition} \label{prop:res}
Let $Z$ be the front projection of a Legendrian link $L$ and let $Y$
be its resolution.  Then $Y$ is the Lagrangian projection of a
Legendrian link $L^\prime$ that is Legendrian isotopic to $L$.
\end{proposition}

\begin{figure}[ht]
\centerline{
\includegraphics[width=2in]{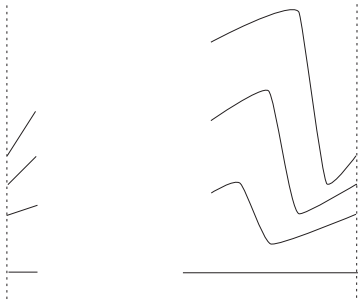}
}
\caption{
Perturbing a front.  The rightmost portion of the figure, which is
included to allow the strands on the right to connect smoothly with
the strands on the left, has Lagrangian projection $T_n$.
}
\label{fig:res-proof}
\end{figure}

\begin{proof}
The corresponding proof in $\R^3$ from \cite{bib:Ng1} distorts the
front, starting at the left and moving right, so that the Lagrangian
projection of the associated Legendrian link is the resolution of
the front. This proof works here until we reach the right end of the
front.  In order to connect smoothly with the left end of the front,
we need to move the right end so that it will close up with the left
end under the identification.  If we label the strands on the right
end of the front by $1,\ldots,n$ in order of increasing $z$
coordinate, then we interpolate a segment of nonpositive slope into
strand $1$ so that it will close up under the identification; next
we interpolate a segment of more negative slope into strand $2$ so
that it will also close up; and so on. See
Figure~\ref{fig:res-proof}. The resulting Lagrangian projection
contains the tangle $T_n$ on its right, and is precisely the
resolution of the original front.
\end{proof}

\begin{figure}[ht]
\centerline{
\includegraphics[width=3in]{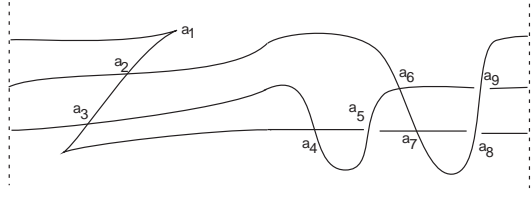}
}
\caption{
The resolved front projection for the front from
Figure~\ref{fig:res-ex} with vertices labeled.  The vertices $a_5$,
$a_8$, and $a_9$ are special crossings.}
\label{fig:dga-ex}
\end{figure}

To calculate the DGA of a Legendrian link $L$ in $\jetS$, it will be
convenient to use a diagram that is a combination of the front and
Lagrangian projections.  First we set some notation. Consider the
projection obtained from the front projection $Z$ of $L$ by formally
appending on the right the tangle $T_n$ for appropriate $n$. We will
call this the {\it resolved front projection} and will often denote
it by $Z^\prime$. See Figure~\ref{fig:dga-ex}. Define the {\it
vertices} of $Z'$ to be its crossings (including the crossings in
$T_n$) and right cusps, and label the vertices by $a_1,\ldots,a_k$.
The term {\it normal crossing} will be used for vertices in the
front projection and for crossings in $T_n$ that appear with the
overstrand having lesser slope. A {\it special crossing} will denote
a crossing in $T_n$ where the overstrand has greater slope.  In the
resolved front projection, all special crossings will be drawn with
a broken segment to denote the understrand, to emphasize the fact
that these are the only crossings which are not resolved by having
the strand with lesser slope be the overstrand.

The DGA for $L$ is $(\A,\d)$, where $\A$ is the free,
noncommutative, unital algebra generated by $a_1,\ldots,a_k$, with
grading to be described in Section~\ref{ssec:gradings}, and $\d$
is defined over $\Z/2$ by the following procedure.
As in \cite{bib:Ch} and
\cite{bib:Ng1}, we define $\d(a_i)$ for a generator $a_i$ by
considering a certain class of immersed disks in the resolved front
projection $Z^\prime$.

\begin{figure}[ht]
\centerline{
\includegraphics[width=4in]{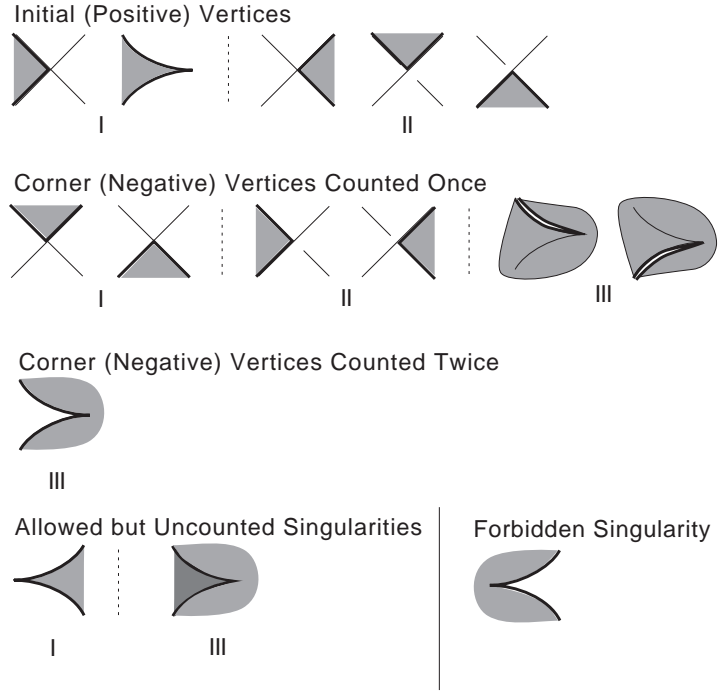}
}
\caption{
Possible singularities in an admissible map and their
classification. The shaded area is the image of the map restricted
to a neighborhood of the singularity;  the heavy line indicates the
image of the boundary of $D^2$.  In two of the diagrams, the heavy
line has been shifted off of itself for clarity.  The diagram with
heavy shading indicates that the image overlaps itself. The regions
denoted by II are new singularities that do not appear in the theory
of Legendrian knots in $\R^3$; they exist due to the presence of the
tangle $T_n$.  For a simple front, the regions denoted by III do not
occur and can be omitted.
}
\label{fig:sing}
\end{figure}

\begin{definition}  An \textit{admissible map} in $Z^\prime$ is an
immersion from
the $2$-disk $D^2$ to $\R^2$ which maps the boundary of $D^2$ into $Z^\prime$,
and which satisfies the following properties:  the image of the
map near any singularity looks locally like
one of the diagrams in Figure \ref{fig:sing}, and, in the notation of
Figure \ref{fig:sing}, there is
precisely one initial vertex.
\end{definition}

Note that Figure~\ref{fig:sing} is quite similar to Figure~5 in
\cite{bib:Ng1} except that due to the special vertices in the tangle
$T_n$, there are additional possibilities for initial vertices and
corner vertices. Note that the ``forbidden singularity" was also
forbidden in \cite{bib:Ng1}. Here it is forbidden since if some disk
involves this singularity, then it has a rightmost negative corner
(at one of the special crossings),  but one can check that any disk
with a rightmost negative corner at a special crossing has to remain
within the tangle $T_n$ and thus cannot ``escape" the tangle on the
left to reach the cusp singularity.

For each admissible map with initial vertex at $a_i$, we can read off
the singularities of the boundary of the disk counterclockwise, beginning
just after $a_i$ and ending just before reaching $a_i$ again. Count each
of the singularities $0$, $1$, or $2$ times, depending on which singularity
it depicts in Figure~\ref{fig:sing}. Concatenating the labels for the
singularities, with multiplicity, yields a word in $\A$, which we associate
to the disk. We then define $\d a_i$ to be the sum of these words over
all admissible maps with initial vertex at $a_i$, or $1$ plus this sum if
$a_i$ is a right cusp.

As is usual, we can extend $\d$ to all of $\A$ via the Leibniz rule.
It is straightforward to check that the usual definition of the DGA
for Lagrangian projections translates to our definition for fronts.
Theorems~\ref{thm:d2} and~\ref{thm:inv} then imply that $\d^2=0$ and
that the DGA $(\A,\d)$ is, up to equivalence, a Legendrian-isotopy
invariant of $L$.

For example, for the resolved front in Figure~\ref{fig:dga-ex},
we have
\begin{alignat*}{3}
\d a_1 &= 1 + a_9a_2 + a_8a_3 & \quad\quad
\d a_4 &= a_3+a_5 & \quad\quad
\d a_7 &= a_6a_5 + a_2 + a_5a_4 + a_8 \\
\d a_2 &= a_5a_3 &
\d a_5 &= 0 &
\d a_8 &= a_9a_5 \\
\d a_3 &= 0 &
\d a_6 &= a_5+a_9 &
\d a_9 &= 0.
\end{alignat*}

\begin{remark}
As in \cite{bib:ENS,bib:Ng1}, we may also
lift the DGA to an algebra over a ring of the form
$\Z[t_1,t_1^{-1},\ldots,t_m,t_m^{-1}]$, graded over $\Z$.
Since we currently have no applications of this lifted DGA, we
will omit its definition here.
\end{remark}

At times, it will be convenient to restrict our attention to a
special class of fronts; as in \cite{bib:Ng1}, it is easy to
Legendrian-isotop any front into the form described below.

\begin{definition} \label{def:simple}
Any front in $\jetS$ that represents a nontrivial element of $\pi_1(\jetS)$
divides the cylinder $S^1
\times
\R$ into several regions, including two unbounded ones, one above,
one below the
front.  We call a front {\it simple} if each of its right cusps
lies on the boundary of one of these unbounded regions.
\end{definition}

The front in Figure~\ref{fig:res-ex} is simple, while the front in
Figure~\ref{fig:stand_calc} below is not. For simple fronts, the
vertices labeled with III in Figure~\ref{fig:sing} do not occur and
thus it follows that the boundary map for the generators $a_i$ can
be calculated in terms of embedded (rather than immersed) disks in
the resolved front diagram.

\subsection{Gradings}
\label{ssec:gradings}

In the previous section, we defined the Chekanov--Eliashberg DGA
associated to a Legendrian front in $\jetS$, except for the grading
on the algebra; we now address the issue of gradings. In
\cite{bib:Ch}, Chekanov defined a grading for generators of the DGA
of a link in $\R^3$.  Following \cite{bib:Ng1}, we will describe how
to define the degrees of vertices in a resolved front projection in
$\jetS$.

For knots, the grading is well-defined, but for links, it depends
on the choice of some auxiliary points on each component of the
link; changing one of the points shifts some of the indices by some
integer. The cases of interest to us will be knots and particular
links for which there is a canonical way to choose auxiliary
points. We address these cases in turn; the definition of gradings
for a general link, which we omit here, can be inferred from the
corresponding definition for links in $\R^3$ in \cite{bib:Ng1}.

For the front of a Legendrian knot in $\jetS$, we can define the
grading on the DGA as follows. The degree of each right cusp is $1$.
For a crossing $a$ in the resolved front projection, traverse the
knot beginning at the greater-slope strand at the crossing, and
ending when the crossing is reached again, and let $c(a)$ be the
number of cusps traversed upwards minus the number traversed
downwards. If $a$ is a normal crossing, then $\deg a = c(a)$; if $a$
is a special crossing, then $\deg a = -c(a)-1$. By extending degree
multiplicatively (i.e., $\deg(ab) = \deg a + \deg b$), this
determines a $\Z/(2r)$ grading on the DGA, where $r$ is the rotation
number of the knot.

Now suppose instead that we are given a Legendrian link $L =
(\Lambda_1,\dots,\Lambda_n)$ in $\jetS$, such that each component
$\Lambda_i$ is Legendrian isotopic to the $1$-jet (in the front
projection, the graph) of a function. We will see that there is a
canonical way to choose an auxiliary point on each component, and
hence the DGA for $L$ can be given a natural grading which is
unique.

\begin{definition} \label{def:initial}
Given a Legendrian knot $\Lambda \subset \mathcal J^1(S^1)$, the
{\it branches} of $\Lambda$ are the connected components of $\Lambda
\backslash C$, where $C$ denotes the set of points that
front-project to cusps. Suppose further that $\Lambda \subset
\mathcal J^1(S^1)$ is Legendrian isotopic to the $1$-jet of the
$0$-function, $\jetO$.
A branch $\mathcal I$ is called {\it initial} if there exists $v \in
\mathcal I$ and a contact isotopy $\kappa_t$ of $\jetS$, $t \in
[0,1]$, so that $\kappa_0=\textrm{id}$, $\kappa_1(\Lambda) = \jetO$,
and $\kappa_t(v)$ never front-projects to a cusp (including the
singular point of a ``Reidemeister type I'' move in the isotopy
which creates or destroys two cusps and a crossing).
\end{definition}

For each crossing $a_i$, let $N^o(a_i)$ and $N^u(a_i)$
denote neighborhoods of $a_i$ on the two strands intersecting
at $a_i$ with $N^o(a_i)$ being a portion of the overstrand.
Note that for a normal crossing, $N^o(a_i)$ will be
in the strand of lesser slope.
%
%
For each $j= 1, \dots, k$, fix a base point $p_j$ on an initial branch
of $\Lambda_j$, and let $P = \{ p_1, \dots, p_k \}$ denote this set
of base points.

To each crossing, we associate two ``capping paths'' $\gamma_i^o,
\gamma_i^u$: $\gamma_i^o$ is a path from $a_i$ through $N^o(a_i)$
to a point in $P$
while
$\gamma_i^u$ is a path from $a_i$ through $N^u(a_i)$ to a point in $P$.
For each such path $\gamma$, let $c(\gamma)$ denote the number of
up cusps minus the number of down cusps along the path.

\begin{definition}
Suppose $L = (\Lambda_1, \dots, \Lambda_n)$ where $\Lambda_i$ is
Legendrian isotopic to the $1$-jet of a function.  Then
$$\deg(a_i) =
\begin{cases}
     1, &a_i \text{ a right cusp}\\
    c(\gamma_i^u) - c(\gamma_i^o), & a_i \text{ a normal crossing}\\
c(\gamma_i^u) - c(\gamma_i^o) - 1,  & a_i \text{ a special
crossing.}
\end{cases}
$$
As usual, degree extends to the entire algebra by
multiplicativity.
\end{definition}

\begin{lemma}  The degree function does not depend
on the choice of paths $\gamma_i^o$, $\gamma_i^u$, nor on the
choice of marked points on initial branches.
\end{lemma}

\begin{proof} Since each component $\Lambda_i$ of the link is
Legendrian isotopic to $j^1(0)$, each component must contain equal
numbers of up and down cusps. It is then easy to check that
$\deg(a_i)$ does not depend on the choice of paths $\gamma_i^o$,
$\gamma_i^u$. In addition, Proposition 5.3 in \cite{bib:Tr}  shows
that a path between any two points in initial branches contains an
equal number of up and down cusps.  Thus the degree function does
not depend on the choice of points in initial branches.
\end{proof}

Notice that for any Legendrian isotopy of a link $L = (\Lambda_1,
\dots, \Lambda_k)$, a choice of base points on initial branches for
the original link yields a choice of base points on initial branches
throughout the isotopy. It follows that the degree given above makes
the DGA of the link a Legendrian-isotopy invariant.

For both knots and the links we have considered in this section,
we have the following result, whose proof can (essentially) be
found in \cite{bib:Ch}.

\begin{proposition} $\partial$ lowers degree by $1$.
\end{proposition}

\subsection{Link DGAs}
\label{ssec:lk_mod}

In addition to the grading of the differential algebra, there is a
structure on the DGA for a link $L$ given by Mishachev's relative
homotopy splitting \cite{bib:Mi}; see \cite{bib:Ng1} for the
formulation we present here. For a link $L = (\Lambda_1, \dots,
\Lambda_k)$, we can split something which is essentially a submodule
of the DGA into $k^2$ pieces which are invariant under Legendrian
isotopy. This additional structure will be useful in applying the
DGA to distinguish between links.

Recall that for a crossing $a$, $N^o(a)$ and $N^u(a)$ denote
neighborhoods of $a$ in the overstrand and understrand.

\begin{definition} For $L = (\Lambda_1, \dots, \Lambda_k)$, let $j_1,
j_2 \in \{1, \dots, k\}$.
If $j_1 \neq j_2$,
define $\A^{j_1 j_2}$ to be the module over $\Z/2$ generated by words of the
form $a_{i_1} \cdots a_{i_m}$ where $N^o(a_{i_1}) \subset \Lambda_{j_1}$,
    $N^u(a_{i_m}) \subset \Lambda_{j_2}$, and
    for $1\leq p<m$, $N^u(a_{i_p})$ and $N^o(a_{i_{p+1}})$ belong to
    the same component of $L$.
For $j_1 = j_2 = j$, let $\A^{j_1j_2}$ be the module generated by such words,
together with an indeterminate $e_j$.  Then let $\A^{**} = \oplus \A^{j_1j_2}$.
\end{definition}

Each $1$ term in the definition of $\partial$ will be replaced by
an $e_j$ term; see below. Note that a generator $a$ is in $\A^{jk}$
if and only if $N^o(a) \subset \Lambda_j$ and $N^u(a) \subset \Lambda_k$.

There is an algebra structure on $\A^{**}$, where multiplication
is defined by the following map $\A^{j_1j_2} \times \A^{j_3j_4}
\rightarrow \A^{j_1j_4}$: the map is $0$ unless $j_2=j_3$, in
which case it is given on generators by concatenation, with the
$e_j$ terms acting as the identity. There exists a differential on
$\A^{**}$ which is a slight variation of $\partial$. Define
$\partial^\prime a_i$ on the generators $a_i$ of $\A^{**}$ as
follows: if $N^o(a_i)$ and $N^u(a_i)$ are contained in distinct
components of $L$, then $\partial^\prime a_i = \partial a_i$; if
$N^o(a_i)$ and $N^u(a_i)$ are contained in the same component
$\Lambda_j$ of $L$, replace any $1$ term by $e_j$. It is easy to
see (cf.\ \cite{bib:Ng1}) that $\partial'$ preserves $\A^{j_1j_2}$
for any $j_1,j_2$.

%

\begin{definition} The {\it link DGA} of $L =
(\Lambda_1, \dots, \Lambda_k)$ is $(\A^{**}, \partial^\prime)$,
where $\partial^\prime a_i$ is defined above, and we extend
$\partial^\prime$ to $\A^{**}$ by applying the Leibniz rule and
setting $\partial^\prime e_j = 0$ for all $j$.  $\A^{**}$ inherits
a grading from the DGA of $L$ with $\deg e_j = 0$ for all $j$.
\end{definition}

We may define grading-preserving elementary and tame automorphisms
for link DGAs as for DGAs, with the additional stipulation that
all maps must preserve the link structure by preserving
$\A^{j_1j_2}$ for all $j_1, j_2$. Similarly, we may define an
algebraic stabilization of a link DGA, with the additional
stipulation that the two added generators both belong to the same
$\A^{j_1j_2}$.  As usual, we then define two link DGAs to be
equivalent if they are tamely isomorphic after some number of
algebraic stabilizations.

\begin{proposition}  If $L$ and $L^\prime$ are Legendrian isotopic
oriented links, then the link DGAs for $L$ and $L^\prime$ are
equivalent.
\end{proposition}

As in \cite{bib:Ng1}, one can define a \textit{link characteristic
algebra} from a link DGA, which is also a Legendrian-isotopy
invariant up to equivalence. This is defined to be the quotient of
$\A^{**}$ by the subalgebra generated by the image of $\d'$.

\subsection{Split Poincar\'e--Chekanov polynomials}
\label{ssec:split_polys}

As in \cite{bib:Mi,bib:Ng1}, we can derive a set of first-order
Poincar\'e--Chekanov polynomials from the link DGA.  To do this, we
will first need the notion of an augmentation for a link. Recall
that an augmentation for a knot consists of an algebra map $\eps :
\A \to \Z/2$ where $\A$ is the DGA over $\Z/2$ for the knot, $\eps
\circ \partial = 0$, $\eps(1)=1$, and $\eps$ vanishes for any
element of nonzero degree.

\begin{definition} Suppose that each component $\Lambda_i$ of $L =
(\Lambda_1, \dots, \Lambda_k)$ when
considered as a knot has an augmentation $\eps_i$.  Extend these
augmentations to all vertices $a_i$
of $L$ by setting
$$\eps(a_i) =
\begin{cases}
\eps_j(a_i), & \text{ if } a_i\in\A^{jj} \text{ for some } j,\\
0, & \text{ otherwise.}
\end{cases}
$$
An {\it augmentation} of $L$ is any function $\eps$ obtained in this way.
\end{definition}

An augmentation $\eps$, as usual, gives rise to a first-order
Poincar\'e--Chekanov polynomial $\chi_{\eps}(\lambda)$; this
polynomial splits into $k^2$ polynomials
$\chi_\eps^{j_1j_2}(\lambda)$ corresponding to the pieces in
$\A^{j_1j_2}$. More precisely, in the expression for $\d a_i$,
replace each $a_j$ by $a_j+\eps(a_j)$; the result has no zero-order
term in the $a_j$, and we denote the first-order term by
$\d^1_\eps(a_j)$. Then $\d^1$ is a differential on the graded vector
space $V_*$ over $\Z/2$ generated by the $a_i$, and it preserves the
subspaces $V_*^{j_1j_2}$ generated by the $a_i\in\A^{j_1j_2}$.

\begin{definition} \label{def:ChekPolys}
For each augmentation $\eps$,  let
\[
\beta_{k}^{j_1j_2}(\eps)  = \dim \frac{\ker \partial_\eps^1 :
V_k^{j_1j_2} \to V_{k-1}^{j_1j_2}} {\im \partial_\eps^1 :
V_{k+1}^{j_1j_2} \to V_{k}^{j_1j_2}}.
\]
We define the \textit{split Poincar\'e--Chekanov polynomials} of the link
$L$ with respect to $\eps$ to be
\[
\chi_\eps^{j_1j_2}(\lambda)[L] = \sum_k \beta_k^{j_1j_2}(\eps)
\lambda^k.
\]
\end{definition}

\noindent
In the case when $L$ is a knot, there is just one polynomial per augmentation,
which we call the Poincar\'e--Chekanov polynomial and denote
$\chi_\eps(\lambda)[L]$.

Of course, the split polynomials depend on a choice of grading of
the link, which is ambiguous in general but is well-defined for
the cases which are of interest to us. The importance of the split
Poincar\'e--Chekanov polynomials derives from the following
result.

\begin{proposition}
If $L$ is a knot, then the set of Poincar\'e--Chekanov polynomials
for all possible augmentations of $L$ is invariant under Legendrian
isotopy. If $L$ is a link whose components are all Legendrian
isotopic to $\jetO$, then for any $j_1,j_2$, the set of split
Poincar\'e--Chekanov polynomials $\chi_\eps^{j_1j_2}(\lambda)[L]$
for all possible augmentations $\eps$ of a link $L$ is invariant
under Legendrian isotopy.
\end{proposition}

\begin{proof}
See \cite{bib:Ch,bib:Ng1}; note that a similar result holds for
arbitrary Legendrian links in $\jetS$.
\end{proof}

\noindent
If a knot or link has only one Poincar\'e--Chekanov polynomial,
then we will sometimes suppress the subscript $\eps$ and write
$\chi(\lambda)[L]$ or $\chi^{j_1j_2}(\lambda)[L]$.

We finish this section by noting one obvious property of split
Poincar\'e--Chekanov polynomials, also mentioned in \cite{bib:Ng1},
which we will use later.

    \begin{proposition} \label{prop:PermRelns}   For $L = (\Lambda_1,
\Lambda_2, \dots, \Lambda_k)$,
let $\sigma$ be a permutation of $\{1, \dots, k\}$ and consider
$L_\sigma = (\Lambda_{\sigma(1)}, \dots, \Lambda_{\sigma(k)})$.  Then
if $L$ has a unique Poincar\'e--Chekanov polynomial,
$$\chi^{\sigma(j_1)\sigma(j_2)}(\lambda)[L_\sigma] =
\chi^{j_1j_2}(\lambda)[L].$$
    \end{proposition}


\section{Legendrian Knots in the Solid Torus}
\label{sec:knots}

In this section, we give several examples of Legendrian knots in the
solid torus that can be distinguished using the DGA invariant. Note
that the generating-function results from \cite{bib:Tr}, which apply
only to two-component links, cannot be used here.

\begin{figure}[ht]
\centerline{
\includegraphics[height=4in]{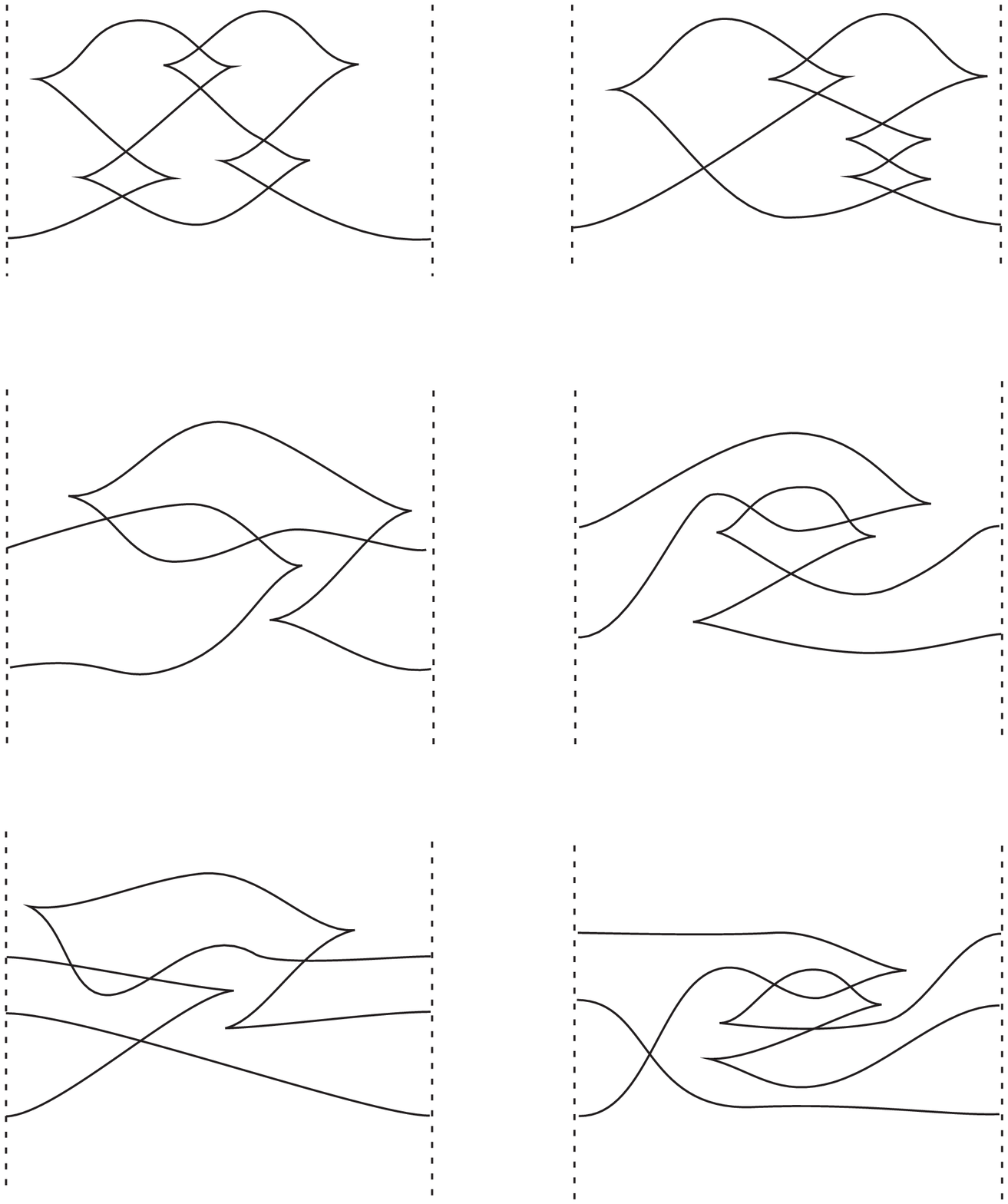}
} \caption{ Pairs of nonisotopic knots in $\jetS$. }
\label{fig:jetspaceknots}
\end{figure}

A large class of Legendrian knots in $\jetS$ can be obtained by
identifying the ends of long Legendrian knots in $\R^3$ (for a
treatment of long Legendrian knots, see \cite{bib:Ch}). More
precisely, given the front of a long Legendrian knot in standard
contact $\R^3$, say with ends asymptotic to the $x$-axis, we can
identify the two ends to produce the front of a Legendrian knot in
$\jetS$; see, for example, the top pair of diagrams in
Figure~\ref{fig:jetspaceknots}, which are derived from long-knot
versions of the Chekanov $5_2$ knots. It is clear from the
resolution construction of Section~\ref{ssec:res} that the DGA for
the solid-torus knot is identical to the DGA of the long knot. One
can then use Poincar\'e--Chekanov polynomials or other techniques to
distinguish these solid-torus knots.

As an example, consider the top pair of knots in
Figure~\ref{fig:jetspaceknots}. One can compute, either directly or
using \cite[Thm.\ 12.4]{bib:Ch} and the standard calculation of the
polynomials for the Legendrian $5_2$ knots, that the knots, and
their long-knot equivalents, have Poincar\'e--Chekanov polynomial
$\chi(\lambda)=2$ (left knot)   and
$\chi(\lambda)=\lambda^2+\lambda^{-2}$ (right knot); hence they are
not Legendrian isotopic.

\begin{figure}[ht]
\centerline{
\includegraphics[width=4.5in]{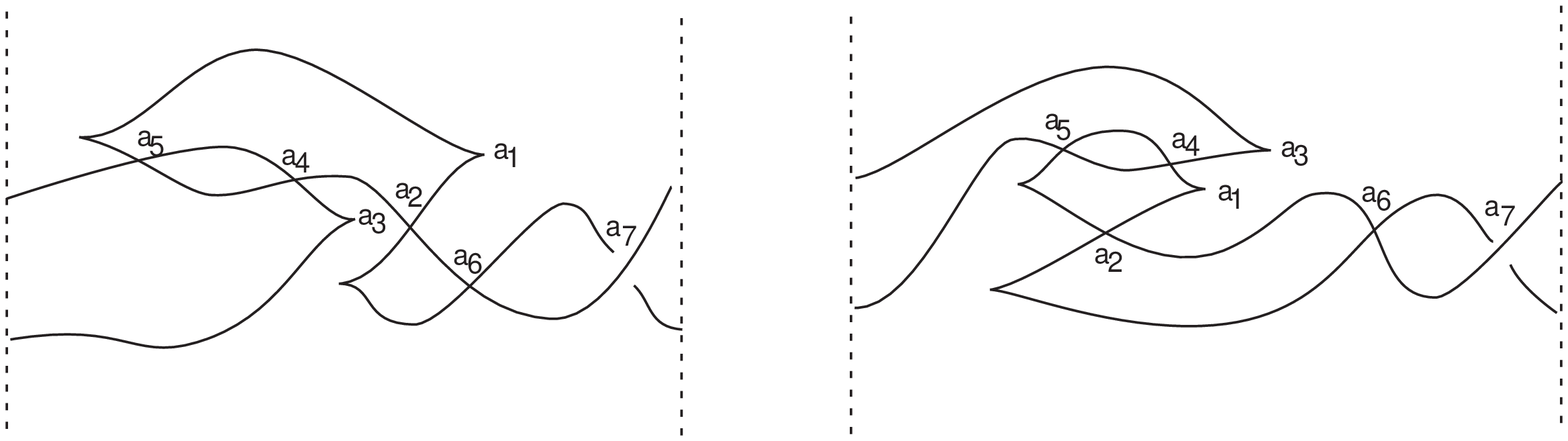}
} \caption{ Calculating the DGA for the middle pair of solid-torus
knots from Figure~\ref{fig:jetspaceknots}. } \label{fig:knot-ex2}
\end{figure}

Another pair of solid-torus knots related to the $5_2$ knots, but
not derived from long knots, is the middle pair of knots in
Figure~\ref{fig:jetspaceknots}, which we will call $K_1$ and
$K_2$. Both $K_1$ and $K_2$ have $r=0$ and $tb=1$. With crossings
as marked in Figure~\ref{fig:knot-ex2}, the DGAs of $K_1$ and
$K_2$ are defined by:
\[
K_1:~~~
\begin{aligned}
\d(a_1) &= 1+(1+a_5a_4)a_2 \\
\d(a_3) &= 1+(1+a_4a_5)a_7 \\
\d(a_6) &= a_2+a_7 \\
\d(a_i) &= 0,~~~i\neq 1,3,6
\end{aligned}
\hspace{0.5in} K_2:~~~
\begin{aligned}
\d(a_1) &= 1+(1+a_4a_5)a_2 \\
\d(a_3) &= 1+a_7(1+a_5a_4) \\
\d(a_6) &= a_2+a_7 \\
\d(a_i) &= 0,~~~i\neq 1,3,6.
\end{aligned}
\]
The key difference between the DGAs for $K_1$ and $K_2$ is that for
$K_1$, $\deg(a_4)=\deg(a_5)=0$, while for $K_2$, $\deg(a_4)=2$ and
$\deg(a_5)=-2$. We can then calculate the Poincar\'e--Chekanov
polynomials to be $\chi(\lambda)[K_1]=\lambda+2$   and
$\chi(\lambda)[K_2]=\lambda^2+\lambda+\lambda^{-2}$. It follows that
$K_1$ and $K_2$ are not Legendrian isotopic.

Using variants of the construction which gives $K_1$ and $K_2$, one
can produce many other examples of pairs of nonisotopic Legendrian
knots. For instance, the bottom pair of knots in
Figure~\ref{fig:jetspaceknots} is a slightly more complicated
example but can also be distinguished using the Poincar\'e--Chekanov
polynomial.


We note that there are other ways to distinguish all of the example
pairs in this section.  One way is the method of Legendrian
satellites (see the Appendix).
A second way to distinguish these pairs is to adapt the generating
function techniques developed for links in \cite{bib:Tr} to the knot
setting; such calculations are done in  \cite{bib:Jo}.

The pair in the second example in Figure~\ref{fig:jetspaceknots} can
be viewed as quotients of Legendrian versions of infinite parity
rational tangles. Quotients of Legendrian versions of even and odd
parity rational tangles are studied in the next section.

\section{Rational Links and Knots}
\label{sec:scp}

In this section, we will calculate the split Poincar\'e--Chekanov
polynomials  of the two-component links studied in \cite{bib:Tr} and
compare these polynomials to the polynomials defined via the theory
of generating functions.  These links can be  seen as quotients of
Legendrian versions of the even parity rational tangles defined by
Conway in \cite{bib:Co}. Some (knot) quotients of Legendrian
versions of odd parity rational tangle will also be examined, thus
leading to solutions of some problems posed in \cite{bib:Tr}.

\begin{figure}[ht]
\centerline{
\includegraphics[height=1.6in]{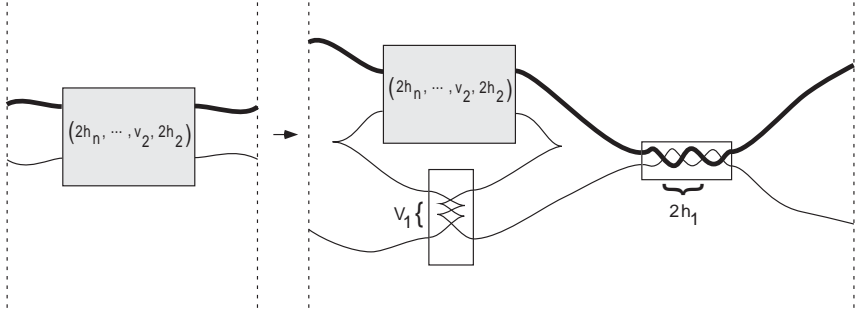}
}
\caption{ The recursive construction of the link
   $(2h_n, v_{n-1}, \dots, 2h_2, v_1, 2h_1)$  from $(2h_n,
v_{n-1}, \dots,
2h_2)$.}
\label{fig:rec_construct}
\end{figure}

The {\it standard rational Legendrian link $(2h_n, v_{n-1}, \dots,
2h_2, v_1, 2h_1)$} is a two-component link that can be constructed
recursively:  for $n =1$, the links $(2h)$  have a front consisting
of the graphs of two functions $f,g : S^1 \to \R$ that intersect
transversally at $2h$ points; for $n \geq 2$, the $(2n-1)$-length
link $(2h_n, v_{n-1}, \dots, 2h_2, v_1, 2h_1)$ is formed from
``vertical and horizontal additions" to the $(2n-3)$-length
   link
$(2h_n,\dots, v_2, 2h_2)$, as shown in
Figure~\ref{fig:rec_construct}. For these two-component links
$(\Lambda_1, \Lambda_2)$, the top/darker strand depicted in
Figure~\ref{fig:rec_construct} will denote the second component. For
all these links, it will be shown that the split
Poincar\'e--Chekanov polynomials are unique  and thus the  set
notation will be dropped.

The polynomials $\chi^{11}(\lambda)[L]$ and
$\chi^{22}(\lambda)[L]$ are the usual Poincar\'e--Chekanov
polynomials for each individual link component.  Since, by
construction,  each component is isotopic to
the
$1$-jet of a function, $\chi^{11}(\lambda)[L] =
\chi^{22}(\lambda)[L] = 0$. However, the polynomials
$\chi^{21}(\lambda)[L]$ and $\chi^{12}(\lambda)[L]$ will carry
useful information.
\begin{definition}  For a Legendrian link $L = (\Lambda_1,
\Lambda_2)$, define
$$
\chi^-(\lambda)[L] = \chi^{21}(\lambda)[L], \qquad
\chi^+(\lambda)[L] = \chi^{12}(\lambda)[L].
$$
\end{definition}

Throughout this section, as in \cite{bib:Tr},
$\overline{p(\lambda)}$ denotes $p(\lambda^{-1})$ for any
polynomial $p$.

\begin{theorem} \label{thm:StandardForm} Consider the Legendrian link
$L = (2h_n, v_{n-1}, \dots,
v_1, 2h_1)$.  Then
\begin{align*}
\chi^+(\lambda)[L] &= h_1 + h_2\lambda^{v_1} + h_3 \lambda^{v_1 + v_2} + \dots
+ h_n \lambda^{v_1 + \dots + v_{n-1}},\\
    \chi^-(\lambda)[L] &=
\begin{cases}
    \overline{\chi^+(\lambda)[L]}, & h_1 \geq 1\\
\lambda^{-1} + 1 + \overline{\chi^+(\lambda)[L]}, & h_1 = 0.
\end{cases}
\end{align*}
\end{theorem}

\begin{proof} First consider the case of $h_1 > 0$.
Draw the link in standard (nonsimple) position as illustrated in
Figure \ref{fig:stand_calc}. This resolved front projection gives rise to a
link DGA with
    \begin{itemize}
\item $2h_1 + \dots + 2h_n$ generators arising from the horizontal
crossings;  these will be denoted by
\begin{align*}
h_{j_i}^\pm, \quad  &j \in \{ 1, \dots, n\}, \  i \in \{1, \dots,
h_j\},\\ \text{\ with\ } &\deg(h_{1_i}^\pm) = 0, \
\deg(h_{j_i}^\pm) = \pm \sum_{k=1}^{j-1} v_k,
\quad
h_{j_i}^+ \in \A^{12}, h_{j_i}^- \in \A^{21};
\end{align*}
\item $v_1 + \dots + v_{n-1}$ generators arising from the vertical
crossings;  these will be denoted by
\begin{align*}
v_{j_i}, \quad  &j \in \{ 1, \dots, n-1\}, \  i \in \{1, \dots,
v_j\},\\ &\text{\ with\ } \deg(v_{j_i} ) = 0,
\quad
v_{j_i} \in \A^{11};
\end{align*}
\item $(v_1 -1) + \dots + (v_{n-1}-1)$ generators arising from the
cusps that occur when $v_i > 1$;  these will be denoted by
\begin{align*}
m_{j_i}, \quad  &j \in \{ 1, \dots, n-1\}, \  i \in \{1, \dots,
v_j-1\},\\
&\text{\ with\ } \deg(m_{j_i} ) = 1,
\quad
m_{j_i} \in \A^{11};
\end{align*}
\item $n-1$ generators arising from the cusps that occur when
connecting the string of $v_i$ crossings with the  string of
$h_{i+1}$ crossings;  these
will be denoted by
\begin{align*}
c_j, \quad j\in \{ 1, \dots, n-1 \}, \quad
\text{\ with\ } \deg(c_{j} ) = 1,
\quad
c_{j} \in \A^{11};
\end{align*}
\item $2$ generators arising from the crossings that occur in the tangle $T_2$;
these will be denoted by
\begin{align*}
t_0^-, t_{-1}^-,
\text{\ with\ } \deg(t_k^- ) = k, \quad t_k^- \in \A^{21}.
\end{align*}
\end{itemize}

\begin{figure}[ht]
\centerline{
\includegraphics[height=1.6in]{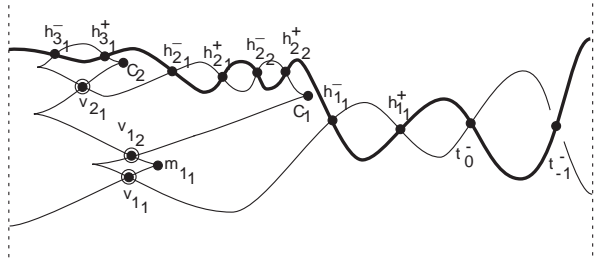}
}
\caption{A nonsimple resolved projection of the link $(2,1,4,2,2)$.}
\label{fig:stand_calc}
\end{figure}

\noindent
Note that to calculate the degrees above, we can set base points
on the two components of $L$ to be the ``leftmost'' points on the fronts,
which certainly lie on initial branches.


There is a unique augmentation for the
link given by $\eps(v_{j_i}) = 1$ for all $j_i$ and $\eps(a) =
0$ for all other vertices.
   The augmented vertices (with $\eps=1$)
are denoted with double circles in Figure~\ref{fig:stand_calc}.
To calculate $\chi^\pm(\lambda)$, it suffices to calculate
$\partial_\eps^1(h_{j_i}^\pm)$ and $\partial_\eps^1(t_{k}^-)$.
    To calculate $\partial_\eps^1$,
we count the number of polygons  constructed by switching branches
at vertices of which at most one is not augmented.
    A direct calculation shows that
\begin{alignat*}{2}
\d_\eps^1(t_{0}^-) &= t_{-1}^-, \quad& \d_\eps^1(t_{-1}^-)&=0,\\
\d_\eps^1(h_{1_{1}}^-) &= t_{-1}^-, \quad&\d_\eps^1(h_{1_1}^+)&=0, \\
\d_\eps^1(h_{j_i}^\pm) &= 0
\quad \text{ for  } j_i \neq 1_{1}.
\end{alignat*}
    Notice that when $v_{k} = 1$,
    $\partial_\eps^1(h_{{k}_{1}}^-)$ will contain a summand of the form
$\sum_i 2h_{(k+1)_i}^- \equiv 0$ which comes from counting disks
with the ``allowed but uncounted singularities" labeled by III in
Figure~\ref{fig:sing}.
    The stated calculations of $\chi^\pm(\lambda)[L]$ follow.

When $h_1 = 0$, the only difference is that the $h_{1_i}^\pm$ no
longer exist and $\partial_\eps^1(t_{0}^-) = 2 t_{-1}^- = 0$. The
desired expressions for $\chi^\pm(\lambda)[L]$ follow for this case.
\end{proof}

\begin{remark} \label{rmk:holo_vs_gf} In \cite{bib:Tr}, polynomials denoted
by $\Gamma^\pm$ for the above Legendrian links were constructed via
the technique of generating functions.  For the Legendrian link   $L
= (2h_n,v_{n-1}, \dots, v_1,2h_1)$,
\begin{align*}
\Gamma^-(\lambda)\left[  L\right] &=
h_1 + h_2\lambda^{-v_1} + h_3\lambda^{-v_1-v_2} + \dots +
h_n\lambda^{-v_1-v_2-\dots-v_{n-1}},
\\
\Gamma^+(\lambda)\left[  L \right] &=
\begin{cases}  \lambda \cdot \Gamma^-(\lambda)\left[
    L \right], & h_1 \geq 1 \\
    (1 + \lambda)  +  \lambda \cdot \Gamma^-(\lambda)\left[ L \right],
& h_1 = 0.
\end{cases}
\end{align*}
Thus for these links, $\chi^+(\lambda)[L] = \overline{\Gamma^-(\lambda)[L]}$,
and, if $h_1 > 0$, $\chi^-(\lambda)[L] = \Gamma^-(\lambda)[L]$. The
same relationship
between $\chi^\pm(\lambda)[L]$ and $\Gamma^\pm(\lambda)[L]$ will hold for
``flyped" versions of these standard links discussed below.
\end{remark}

All these standard rational Legendrian links are topologically
unordered. Notice that Proposition~\ref{prop:PermRelns} implies that
if $L = (\Lambda_1, \Lambda_2)$ and
    $\overline{L} = (\Lambda_2, \Lambda_1)$, then
$$ \chi^-(\lambda)[\overline L] = \chi^+(\lambda)[L] \quad
\text{ and } \quad \chi^+(\lambda)[\overline L] = \chi^-(\lambda)[L].$$
    Thus by the formulas in Theorem~\ref{thm:StandardForm}  we
reproduce the following result from \cite{bib:Tr}.

\begin{corollary} Consider $L = (2h_n,v_{n-1}, \dots, v_1,2h_1)$. Then $L =
(\Lambda_1, \Lambda_2)$ is isotopic to $\overline L = (\Lambda_2,
\Lambda_1)$ if and only if $L = (2h_1)$, $h_1 > 0$.
\label{cor:swap}
\end{corollary}

\begin{figure}[ht]
\centerline{
\includegraphics[height=1.6in]{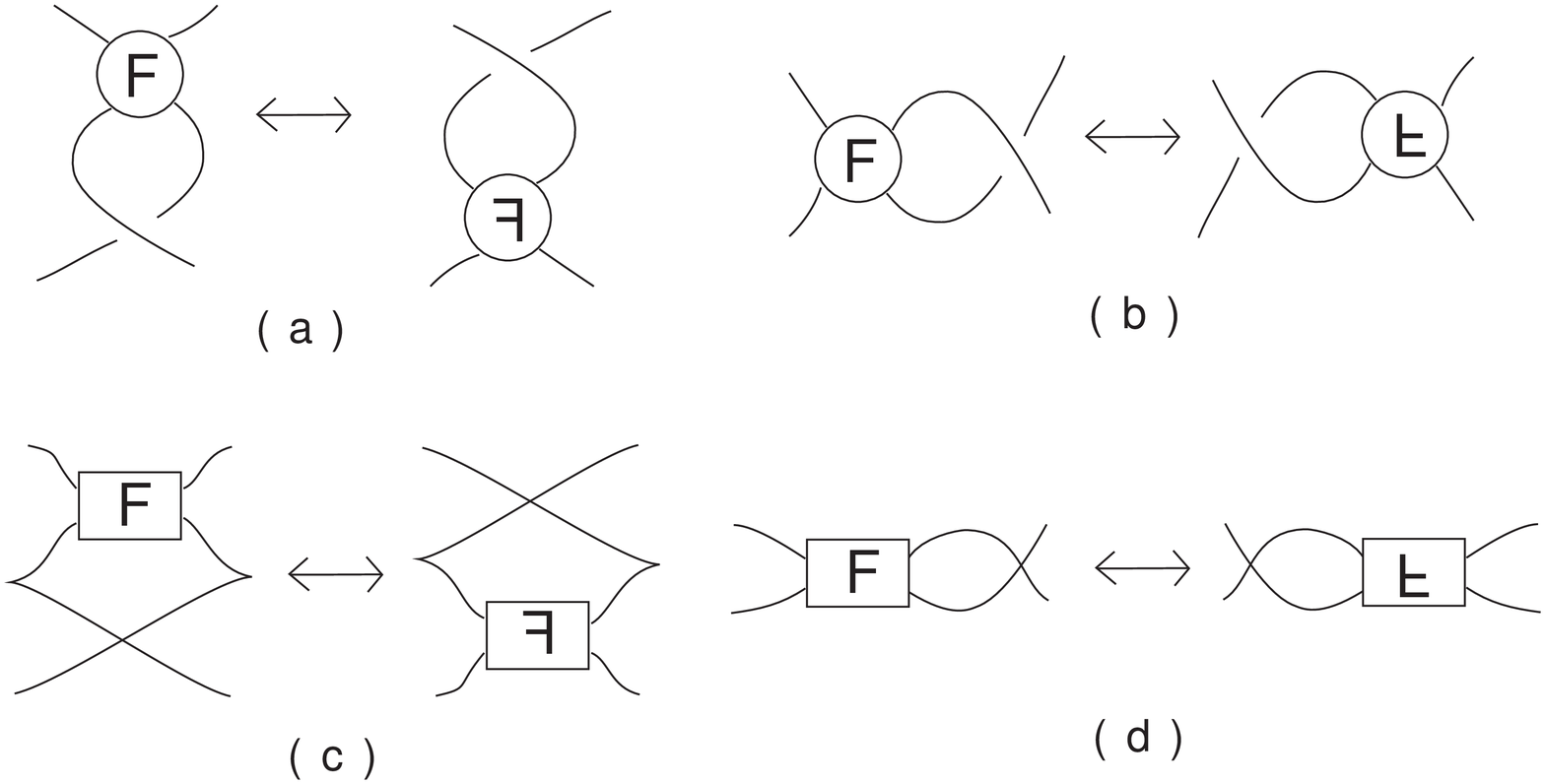}
}
\caption{ (a) A topological vertical flype, (b) a
topological horizontal flype, (c) a Legendrian vertical flype, (d) a
Legendrian horizontal flype.} \label{fig:flype_def}
\end{figure}

We next consider  rational links that are not in the
``standard form'' $(2h_n,v_{n-1},\dots,v_1,2h_1)$, but rather
involve added ``flypes".
For the topological version  of the link $(2h_n, \dots, v_1, 2h_1)$,
$n \geq 2$,  there
   are ``flyping" moves that do not change the topological type of the link.
A topological flype occurs when  a portion of the link, represented
by the circle labeled with ``F" in Figure~\ref{fig:flype_def}, is
rotated $180^\circ$ about a vertical axis (a vertical flype), or
about a horizontal axis (a horizontal flype).   For background on
such flypes, see, for example, \cite{bib:Ad}. This motivates the
definition of a {Legendrian flype}:  when a crossing is formed by
two   edges emanating from  a Legendrian tangle, represented by the
box labeled with ``F" in Figure~\ref{fig:flype_def}, a Legendrian
vertical (horizontal) flype occurs when the tangle is rotated
$180^\circ$ about vertical (horizontal) axis and the crossing  is
``transferred" to the opposite edges.  This rotation action is not a
Legendrian isotopy; hence, although the resulting Legendrian links
are topologically isotopic, they are potentially not Legendrian
isotopic.


For each positive horizontal entry $2h_i$, $i \neq n$, in the
Legendrian link \newline $(2h_n, \dots, 2h_2, v_1, 2h_1)$,  it is
possible to perform $0$, $1$, \dots, or $2h_i$ successive horizontal
flypes; for each vertical entry $v_i$, it is possible to perform
$0$, $1$, \dots, or $v_i$ successive vertical flypes. In
\cite{bib:Tr} it was shown that vertical flypes of a standard
Legendrian link produce an isotopic Legendrian link.  In the
following, we will focus on horizontal flypes of a standard
Legendrian link. The nomenclature
$$
\left(2h_n, v_{n-1}, 2h_{n-1}^{p_{n-1}}, \dots, 2h_2^{p_2}, v_1,
2h_1^{p_1}\right),
   p_i \in \{ 0, \dots, 2h_i\},
$$
will be used to denote the modification of the
   standard  link
by $p_i$ horizontal flypes in the $i^{th}$ horizontal
component.
With this notation, the standard rational link
is written as
$(2h_n, v_{n-1}, 2h_{n-1}^0, \dots, 2h_2^0, v_1, 2h_1^0)$.
If no superscript is
specified for an entry of the vector, it will be assumed to be $0$.
Figure~\ref{fig:flype_calc} illustrates a resolved front
projection of $(2,1,4^1, 2,2)$.

The following theorem shows that the polynomials of these
flyped versions of the standard link also have a nice pattern: as was
seen for the
calculation of the generating function polynomials \cite{bib:Tr},
   the coefficient of  $v_i$ (in the exponent of $\lambda$)
in the formulas given for the polynomials of the standard rational
Legendrian link by Theorem~\ref{thm:StandardForm} will  become for
the flyped versions either $\pm1$ depending on the parity of the
number of flypes in the $h_1, \dots, h_i$ positions.

\begin{theorem} \label{thm:FlypeForm} Let $L$ be the Legendrian
link $\left(2h_n,v_{n-1}, 2h_{n-1}^{p_{n-1}}, \dots,
v_1,2h_1^{p_1}\right).$
     For $j = 1, \dots, n-1$, let
$\sigma(j) = \sum_{i=1}^j p_i$. Then
\begin{align*}
\chi^+(\lambda) [ L ] &= h_1 + \sum_{i=2}^n h_i
\lambda^{(-1)^{\sigma(1)}v_1 + (-1)^{\sigma(2)} v_2+
\dots + (-1)^{\sigma(i-1)}v_{i-1}}, \\
\chi^-(\lambda) [ L ] &=
\begin{cases}
    \overline{\chi^+(\lambda)\left[ L \right]}, & h_1 \geq 1 \\
(\lambda^{-1} + 1) +  \overline{\chi^+(\lambda)[L]}, & h_1 = 0.
\end{cases}
\end{align*}
\end{theorem}

\begin{proof}   By Theorem~\ref{thm:StandardForm}, the stated formulas hold for
$L$ when $p_i = 0$ for all $i$. For arbitrary $p_{n-1}, \dots,
p_1$, assume the formulas hold for
$$L = \left(2h_n,v_{n-1}, 2h_{n-1}^{p_{n-1}}, \dots, v_1,2h_1^{p_1}\right).$$
Consider $L^\prime$ which
differs from $L$ by one additional horizontal flype at the $k^{th}$ position:
\begin{align*}
L^\prime = \left(2h_n,v_{n-1}, 2h_{n-1}^{w_{n-1}}, \dots,
v_1,2h_1^{w_1}\right), \quad \exists k : w_i &= p_i, \text{ for } i \neq k, \\
\text{
and } w_k &= p_k + 1.
\end{align*}
There is a resolved front diagram $Z^\prime$ for $L^\prime$ with the
same number of vertices as a resolved front diagram $Z$ of
     $L$:  in analogue with the argument in the proof of
Theorem~\ref{thm:StandardForm}, let $h_{j_i}^\pm, v_{j_i},
m_{j_i}, c_j, t_k^-$ denote the vertices of $Z$ corresponding to
the horizontal crossings, vertical crossings, crossings resulting
from $v_i > 1$, connections between strings of horizontal and
vertical crossings, and the crossings in the tangle $T_2$.
Similarly, let ${h_{j_i}^\pm}^\prime, {v_{j_i}}^\prime,
m_{j_i}^\prime, c_j^\prime, {t_k^-}'$ denote the corresponding
vertices of $Z^\prime$. For an example, see
Figures~\ref{fig:stand_calc} and \ref{fig:flype_calc}.
\begin{figure}[b]
\centerline{
\includegraphics[height=1.6in]{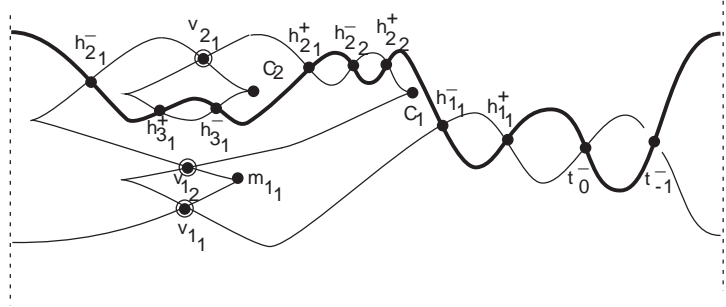}
}
\caption{Labels for a resolved front projection of a flype of the
figure in Figure~\ref{fig:stand_calc}. All vertices are further
labeled by a prime ($'$) which is omitted.}
\label{fig:flype_calc}
\end{figure}
We first calculate the degrees of ${h_{j_*}^\pm}^\prime$ assuming
that the degree of ${h_{j_*}^\pm}$ is $\pm \sum_{\ell=1}^{j-1}
(-1)^{\sigma(\ell)}v_\ell$ where   $\sigma(\ell) =  \sum_{i=1}^\ell
p_i$. It is easy to verify that when $j=1$, $\deg h_{j_*}^\pm = 0$.
When $j>1$, we will show that
     $\deg {h_{j_i}^\pm}^\prime = \pm \sum_{\ell=1}^{j-1}
(-1)^{\sigma^\prime(\ell)}v_\ell$, where
$$\sigma^\prime(\ell) :=  \sum_{j=1}^\ell w_j
=
\begin{cases}
\sigma(\ell), & \ell \leq k-1 \\
\sigma(\ell) + 1, & \ell \geq k.
\end{cases}$$
For each capping path $\gamma$ for $h_{j_i}^\pm$, there is a natural
capping path $\gamma^\prime$ for ${h_{j_i}^\pm}^\prime$ so that $\gamma$ and
$\gamma^\prime$ pass through corresponding sequences of vertical crossings.
    To calculate the degrees,
we must count the number of up and down cusps along these capping
paths.  However, a path passes through a cusp if and only if it
passes through a
vertical crossing.  Thus
it suffices to keep track of how capping paths $\gamma$ and
$\gamma^\prime$ differ in how
they pass through the vertical crossings.
Capping paths $\gamma$ for $h_{j_i}^\pm$ (and thus
$\gamma^\prime$ for ${h_{j_i}^\pm}^\prime$) can be chosen that pass
through each of the vertical vertices $v_{1_{*}}, \dots,
v_{{j-1}_{*}}$ precisely once and do not pass through the vertices
$v_{j_*}, \dots, v_{{n-1}_*}$.  If the flype occurs at the
$k^{th}$ position, cusp counts along paths that do not intersect
$v_{k_*}, v_{{k+1}_*}, \dots, v_{{n-1}_*}$ will be unchanged.
Thus it is easy to see that
$$j \leq k \implies \deg {h_{j_i}^\pm}^\prime = \deg {h_{j_i}^\pm}
= \pm \sum_{\ell=1}^{j-1} (-1)^{\sigma(\ell)}v_\ell = \pm
\sum_{\ell=1}^{j-1} (-1)^{\sigma^\prime(\ell)}v_\ell.$$ A capping
path $\gamma$ for $h_{j_i}^\pm$, $j \geq k+1$, can be broken into
two paths $\gamma_1$ and $\gamma_2$ where $\gamma_1$ passes through
each of the vertical vertices $v_{1_*}, \dots, v_{{k-1}_*}$
precisely once and $\gamma_2$ passes through each of the vertical
vertices $v_{k_*}, \dots, v_{{j-1}_*}$ precisely once. If
$\gamma_1^\prime$ and $\gamma_2^\prime$ denote the corresponding
capping paths for ${h_{j_i}^\pm}^\prime$, then the cusp counts for
$\gamma_1$ and $\gamma_1^\prime$ agree while the cusp counts for
$\gamma_2$ and $\gamma_2^\prime$ differ in sign. Thus if
$$\deg {h_{j_i}^\pm} = \pm \sum_{\ell=1}^{j-1} (-1)^{\sigma(\ell)}v_\ell$$
then
\begin{align*}
\deg {h_{j_i}^\pm}^\prime
&= \deg \gamma_1 - \deg \gamma_2 \\
&= \pm\left(\sum_{\ell=1}^{k-1} (-1)^{\sigma(\ell)}v_\ell +
\sum_{\ell=k}^{j-1} (-1)^{\sigma(\ell)+1}v_\ell \right)= \pm \left(
\sum_{\ell=1}^j
(-1)^{\sigma^\prime(\ell)} v_\ell\right). \\
\end{align*}

It is easy to verify that there are unique augmentations $\eps$ of
$L$ and $\eps^\prime$
of $L^\prime$ given by augmenting precisely the vertical vertices.
We next show that
the calculation of $\partial_\eps^1$ for $L$   will be
``preserved" under the flyping move.
In particular, we will show that when $h_1 > 0$,
\[
\begin{aligned}
\d_\eps^1(t_{0}^-) &= t_{-1}^- & \d_\eps^1(t_{-1}^-)&=0\\
\d_\eps^1(h_{1_{1}}^-) &= t_{-1}^- &\d_\eps^1(h_{1_1}^+)&=0 \\
\d_\eps^1(h_{j_i}^\pm) &= 0
   &\text{ for  } j_i \neq 1_{1}
\end{aligned}
    \implies
\begin{aligned}
\d_{\eps^\prime}^1({t_{0}^-}^\prime) &= {t_{-1}^-}^\prime &
\d_{\eps^\prime}^1({t_{-1}^-}^\prime)&=,\\
\d_{\eps^\prime}^1({h_{1_{1}}^-}^\prime) &= {t_{-1}^-}^\prime
&\d_{\eps^\prime}^1({h_{1_1}^+}^\prime)&=0
\\
\d_{\eps^\prime}^1({h_{j_i}^\pm}^\prime) &= 0
   &\text{ for  } j_i \neq 1_{1}.
\end{aligned}
\]
By an analysis of the flyping procedure and the construction of these links,
   we find that when the flype
occurs at the $k^{th}$ position, if there are no disks with a corner
label of $h_{k_*}$ or $h_{(k+1)_*}$ in the calculation of
$\d_\eps^1(a)$ then $\d_{\eps^\prime}^1(a^\prime)$ will be given by
the formula for $\d_\eps^1(a)$ with primes inserted for all the
variables. This verifies the statement for
$\d_{\eps^\prime}^1({t_{0}^-}^\prime)$,
$\d_{\eps^\prime}^1({t_{-1}^-}^\prime)$, and
$\d_{\eps^\prime}^1({h_{j_i}^\pm}^\prime)$ when $j \neq k, k-1$. A
careful analysis shows that
$\d_{\eps^\prime}^1({h_{(k-1)_*}^\pm}^\prime)$ is given by inserting
primes into the formula for $\d_{\eps}^1({h_{(k-1)_*}^\pm})$, and
that when $v_k > 1$, there are no admissible maps for the
calculation $\d_{\eps }^1({h_{k_*}^\pm})$ coming from disks with
corner labels $h_{k_*}$,  $h_{(k+1)_*}$. When $v_k = 1$, $\d_{\eps
}^1({h_{k_j}^\pm})$ containing a summand of the form $\sum 2
h_{{k+1}_i}^\pm \equiv 0$ implies that
$\d_{\eps^\prime}^1({h_{k_j}^\mp}^\prime)$ contains a summand of the
form $\sum 2 {h_{{k+1}_i}^\mp}^\prime \equiv 0$. (Compare
$\d_\eps^1(h_{2_1}^-)$ in Figure~\ref{fig:stand_calc}, which has two
$h_{3_1}^-$ terms, with $\d_{\eps^\prime}^1({h_{2_1}^+}^\prime)$ in
Figure~\ref{fig:flype_calc}, which has two ${h_{3_1}^+}^\prime$
terms.)
This demonstrates the
claimed calculations of $\d_{\eps^\prime}^1$ when $h_1 > 0$. When
$h_1 = 0$,    $h_{1_i}^\pm$, ${h_{1_i}^-}^\prime$ no longer exist
and the above arguments show that
$$\partial_\eps^1(t_{0}^-) = 2 t_{-1}^- = 0 \implies
\partial_{\eps^\prime}^1({t_{0}^-}^\prime) = 2 {t_{-1}^-}^\prime = 0.$$
  From this, we can
deduce the claimed calculations of $\chi^\pm(\lambda)[L^\prime]$ for
the unique augmentation $\eps^\prime$.
\end{proof}

Theorem~\ref{thm:StandardForm} and Theorem~\ref{thm:FlypeForm} can
be combined to construct numerous examples of nonisotopic Legendrian
links. For example, the links $L = (2,1,2)$ and $L^\prime =
(2,1,2^1)$ shown in Figure~\ref{fig:link+knot_flypes} are not
Legendrian isotopic since
\begin{alignat*}{2}
\chi^+(\lambda)[L] &= 1 + \lambda, &\quad \chi^-(\lambda)[L] & = 1 +
\lambda^{-1}, \\
   \chi^+(\lambda)[L^\prime] &= 1 + \lambda^{-1}, &\quad
\chi^-(\lambda)[L^\prime] & = 1 + \lambda.
\end{alignat*}

If we disregard the splitting of the DGA into modules $\A^{j_1j_2}$,
the standard Poincar\'e--Chekanov polynomial of $L =
\left(2h_n,v_{n-1}, 2h_{n-1}^{p_{n-1}}, \dots,
v_1,2h_1^{p_1}\right)$ is
$$\chi(\lambda)[L] = \chi^-(\lambda)[L] + \chi^+(\lambda)[L].$$
Notice that this nonsplit polynomial would not distinguish between
the links $L$ and $L^\prime$ listed above. However, the nonsplit
polynomial will be useful when examining knots.

Using the above procedure, a vector of the form
$$ \left(2h_n,v_{n-1}, 2h_{n-1}^{p_{n-1}}, \dots,
v_1,(2h_1-1)^{p_1} \right), \quad h_i, v_i \geq 1$$ is a knot that
can be viewed as the quotient of a Legendrian version of an odd
parity rational tangle.

\begin{theorem} \label{thm:knotpoly}
The Legendrian knot $K = \left(2h_n,v_{n-1},  \dots, v_1,(2h_1-1)
\right)$, $h_i, v_i \geq 1$ has unique first-order
Poincar\'e--Chekanov polynomial
$$\chi(\lambda)[K]
    =  \sum_{i=2}^n h_i \lambda^{-v_1 - v_2 -
\dots -  v_{i-1}} + (2h_1 - 1) + \sum_{i=2}^n h_i \lambda^{v_1 +
v_2+ \dots +  v_{i-1}}.$$ The Legendrian knot with flypes $\widehat
K = \left(2h_n,v_{n-1}, 2h_{n-1}^{p_{n-1}}, \dots,
v_1,(2h_1-1)^{p_1}\right)$ has unique first-order
Poincar\'e--Chekanov polynomial
\begin{align*}
\chi(\lambda)[\widehat K] &= \sum_{i=2}^n h_i
\lambda^{-((-1)^{\sigma(1)}v_1 + (-1)^{\sigma(2)} v_2+
\dots + (-1)^{\sigma(i-1)}v_{i-1})}   + (2h_1-1) \\
&\quad + \sum_{i=2}^n h_i \lambda^{(-1)^{\sigma(1)}v_1 +
(-1)^{\sigma(2)} v_2+ \dots + (-1)^{\sigma(i-1)}v_{i-1}},
\end{align*}
where $\sigma(j) = \sum_{i=1}^j p_i$ for $j=1,\dots,n-1$.
\end{theorem}

This theorem is proven with calculations analogous to those in the
proofs of Theorems~\ref{thm:StandardForm} and  \ref{thm:FlypeForm}.
Label crossings similarly to the two-component case, but without the
$\pm$ superscripts. Essentially, one uses the augmentation which
augments precisely the $v_{j_i}$ generators, and computes the
linearized differential for the $h_{j_i}$ and $t_i$ generators as
before. We now need to consider the linearized differential on the
$c_j$, $m_{j_i}$, and $v_{j_i}$ generators as well; the fact that
these make no contribution to linearized homology follows from the
analogous fact, in the two-component link setting, that the split
Poincar\'e--Chekanov polynomial $\chi^{11}(\lambda)$ is $0$.

The main difficulty in the knot case is that there may be other
augmentations besides the one used above. For example, in the
standard form, any augmentation must augment the vertical vertices
but it is now possible to augment any among the set $\{ h_{1_*},
t_0\}$. (Recall that for the link situation, these crossings were
interstrand crossings and thus they could not be augmented.)  It is
easy to verify that changing the augmentations of $\{ h_{1_*},
t_0\}$ will not change the  admissible maps used in the calculation
of $\d_\eps^1$ and thus $\chi_\eps(\lambda)[K]$ will be independent
of the choice of augmentation. For the flyped knots $\widehat K$, it
is possible to augment some generators ${h_{j_*}}$ when
$\deg({h_{j_*}}) = 0$. This may result in an addition of immersed
disks in the calculation of $\d_{\eps}^1$, but such disks always
occur in canceling pairs and thus the polynomials
$\chi_\eps(\lambda)[\widehat K]$ are independent of augmentation.

We can apply Theorem~\ref{thm:knotpoly} to answer Question 1.34 from
\cite{bib:Tr} in the negative: the Legendrian knots $(2,1,2,1,1)$
and $(2,1,2^1,1,1)$ (see Figure~\ref{fig:link+knot_flypes}) are not
Legendrian isotopic, because they have Poincar\'e--Chekanov
polynomial $\lambda^2+\lambda+1+\lambda^{-1}+\lambda^{-2}$ and
$\lambda+3+\lambda^{-1}$, respectively.

\section{Connect Sums}
\label{sec:connectsums}

\begin{figure}[ht]
\centerline{
\includegraphics[height=1in]{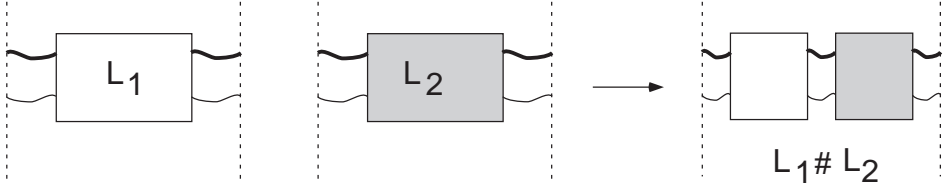}
}
\caption{ The construction of the connect sum $L_1 \# L_2$.}
\label{fig:connect_sum}
\end{figure}

In this section, we consider (usually nonrational) ``connect sums" of
rational knots and links,
in the sense introduced in \cite{bib:Tr}.
     Since, up to Legendrian
isotopy, the connect sum may depend on the choice of where
the links are cut into tangles, a standard position for
   cutting the links will be chosen.  Namely, the connect sum
$L_1 \# L_2$  is defined as the closure of the connect sum of the
Legendrian rational tangles  $\mathcal T(L_1)$ and $\mathcal
T(L_2)$, which are constructed analogously to the links $L_i$.  This
construction is illustrated in Figure~\ref{fig:connect_sum} where,
if $L_1$ denotes the link $(2h_n, \dots, 2h_1)$, then $\mathcal
T(L_1)$ corresponds to  Figure~\ref{fig:rec_construct} except
considered as a tangle rather than closed to a link. When working
with the polynomials constructed from the theory of generating
functions, simple formulas exist for constructing the polynomials of
$L_1 \# L_2$ from the polynomials of $L_1$ and $L_2$.  Similarly,
simple formulas hold for the polynomials constructed from the
holomorphic theory.

\begin{theorem}  \label{thm:connect_sum}
Consider the Legendrian links
$$
    L_1 = (2h_n, v_{n-1}, 2h_{n-1}^{p_{n-1}}, \dots, v_1,
    2h_1^{p_1}),~
    L_2 = (2k_m, u_{m-1}, 2k_{m-1}^{w_{m-1}}, \dots, u_1, 2k_1^{w_1}).
$$
    Then
\begin{align*}
\chi^+(\lambda)[L_1 \# L_2] &= \chi^+(\lambda)[L_1] + \chi^+(\lambda)[L_2];\\
\chi^-(\lambda)[L_1 \# L_2] &=
\begin{cases}
    \chi^-(\lambda)[L_1] + \chi^-(\lambda)[L_2] , & h_1, k_1 \geq 1 \\
    \chi^-(\lambda)[L_1] + \chi^-(\lambda)[L_2]  - (\lambda^{-1} + 1), &
\text{else.}
\end{cases}
\end{align*}
\end{theorem}

\begin{proof}
Let $Z_1^\prime, Z_2^\prime, Z_{1+2}^\prime$ denote resolved front
projections for $L_1, L_2,L_1 \# L_2$, and then let $V^{12}_1,
V^{12}_2, V^{12}_{1+2}$ and $V^{21}_1, V^{21}_2, V^{21}_{1+2}$
denote the first-order summands of the corresponding link DGAs.
Then, as in the proof of Theorem~\ref{thm:FlypeForm}, if $V_1^{12}$
is generated by $h_{j_i}^+$ and $V_2^{12}$ is generated by
$k_{m_\ell}^+$, then $V_{1+2}^{12}$ is generated by $h_{j_i}^+,
k_{m_\ell}^+$; if $V_1^{21}$ is generated by $h_{j_i}^-, t_0^-,
t_{-1}^-$ and $V_2^{21}$ is generated by $k_{m_\ell}^-, u_0^-,
u_{-1}^-$, then $V_{1+2}^{21}$ is generated by $h_{j_i}^-,
k_{m_\ell}^-, u_0^-, u_{-1}^-$.  Notice that the two generators
$t_{0}^-, t_{-1}^-$ that were needed to close up $L_1$ are no longer
necessary due to the presence of $L_2$.

It is easily verified that augmentations $\eps_1$ and $\eps_2$ of
$Z_1^\prime$ and $Z_2^\prime$ lead to an augmentation
$\eps_{1+2}$ of $Z_{1+2}^\prime$ (and conversely) and that on
$V_{1+2}^{12}$, $\partial_{\eps_{1+2}}^1 =
\partial_{\eps_{1}}^1 +  \partial_{\eps_{2}}^1$, where we extend
$\d_{\eps_1}^1,\d_{\eps_2}^1$ trivially to $V_2^{12},V_1^{12}$
respectively. It follows
that $\chi^+(\lambda)[L_1 \# L_2] = \chi^+(\lambda)[L_1] +
\chi^+(\lambda)[L_2]$.

Recall from the proof of Theorem~\ref{thm:FlypeForm} that on
$V_1^{21}$ and $V_2^{21}$,
    \begin{align*}
\ker \partial_{\eps_1}^1/\im \d_{\eps_1}^1 \text{ is spanned by }
&\begin{cases}
h_{1_{1}}^- + t_0^-, h_{j_i}^-, &   j_i \neq 1_{1}, \text{ when }h_1 \geq 1 \\
h_{j_i}^-, t_0^-, t_{-1}^-,  & h_1 = 0;
\end{cases}\\
\ker \partial_{\eps_2}^1/\im \d_{\eps_2}^1 \text{ is spanned by }
&\begin{cases}
k_{1_{1}}^- + u_0^-, k_{m_\ell}^-, &   m_\ell \neq 1_{1}, \text{ when
}k_1 \geq 1\\
k_{m_\ell}^-, u_0^-, u_{-1}^-,  &  k_1 = 0. \\
\end{cases}
\end{align*}
If $h_1= 0$ or $k_1 = 0$, we may assume by applying a translation
if necessary that $h_1 = 0$.  Then we see that
$\partial_{\eps_{1+2}}^1 = \partial_{\eps_1}^1 + \partial_{\eps_2}^1$ on
the generators of $V_{1+2}^{21}$ and thus
$\chi^-(\lambda)[L_1 \# L_2]$ equals
$\chi^-(\lambda)[L_1] + \chi^-(\lambda)[L_2]$ minus two terms of
degree $-1$ and $0$.
Else if $h_1, k_1 > 0$, then it is easy to verify that on  $V_{1+2}^{21}$,
$\ker \partial_{\eps_{1+2}}^1/\im \d_{\eps_{1+2}}^1$  is spanned by
$$
h_{1_{1}}^- + u_0^-, h_{j_i}^-, k_{m_\ell}^-, \quad   j_i \neq 1_{1};
$$
it follows that $\chi^-(\lambda)[L_1 \# L_2] = \chi^-(\lambda)[L_1] +
\chi^-(\lambda)[L_2]$.
    \end{proof}

\begin{remark} \label{rmk:gfconnectsum}
   An analogous formula for connect sums holds for the polynomials
constructed by the theory of generating functions.  Namely, for the
Legendrian links considered in the statement of
Theorem~\ref{thm:connect_sum},
we have
\begin{align*}
\Gamma^-(\lambda)[L_1 \# L_2] &= \Gamma^-(\lambda)[L_1] +
\Gamma^-(\lambda)[L_2];\\
\Gamma^+(\lambda)[L_1 \# L_2] &=
\begin{cases}
    \Gamma^+(\lambda)[L_1] + \Gamma^+(\lambda)[L_2] , & h_1, k_1 \geq 1 \\
    \Gamma^+(\lambda)[L_1] + \Gamma^+(\lambda)[L_2]  - (1 +\lambda), &
\text{ else.}
\end{cases}
\end{align*}
Thus, in parallel to what we observed for  rational links,
\begin{align*}
\chi^+(\lambda)[L_1 \# L_2] &= \overline{\Gamma^-(\lambda)[L_1 \# L_2]},\\
\chi^-(\lambda)[L_1 \# L_2] &= \Gamma^-(\lambda)[L_1 \# L_2], \qquad
\text{ when } h_1, k_1 \geq 1.
\end{align*}
\end{remark}

We can use the same techniques to calculate the Poincar\'e--Chekanov
polynomial for the connect sum of a rational knot and a
rational link.

\begin{theorem} \label{thm:connectsumknot} Consider the Legendrian
link and knot
\begin{align*}
    L_1 &= (2h_n, v_{n-1}, 2h_{n-1}^{p_{n-1}}, \dots, v_1, 2h_1^{p_1}),\\
    L_2 &= (2k_m, u_{m-1}, 2k_{m-1}^{w_{m-1}}, \dots, u_1, (2k_1-1)^{w_1}).
\end{align*}
    Then
$$
\chi(\lambda)[L_1 \# L_2] =
\begin{cases}
\chi^+(\lambda)[L_1] + \chi^-(\lambda)[L_1]
+\chi(\lambda)[L_2], & h_1 \geq 1\\
\chi^+(\lambda)[L_1] + \chi^-(\lambda)[L_1] - (\lambda^{-1} + 1)
+\chi(\lambda)[L_2], & h_1 = 0.
\end{cases}
$$
\end{theorem}


Next we will explore some examples of connect sums that can and
cannot be distinguished by their Poincar\'e--Chekanov polynomials.
Then it will be shown that sometimes the characteristic algebra can
distinguish links that cannot be differentiated by their
polynomials. For simplicity, we will restrict to connect sums made
from ``basic building blocks''.

\begin{definition}
Let $\mathcal L_{j_1, j_2}$ denote
the connect sum $(2,1,j_1) \# (2,1,j_2)$ as shown in
Figure~\ref{fig:L_{j_1,j_2}}.
\end{definition}

Notice that there
are many ways to label a connect sum.  For example,
$\mathcal L_{2, 2} = (2,1,2) \# (2,1,2)$ can also be written
as
$(2,1,3^1) \# (2,1,1)$ or as $(2,1,4^2) \# (2,1,0)$.
The notation $\mathcal L_{j_1, j_2}$ is convenient
since it emphasizes the geometrical aspect that the
link is made out of two basic ``bubble units" separated by
$j_1$ and $j_2$ crossings.

\begin{figure}[ht]
\centerline{
\includegraphics[height=1in]{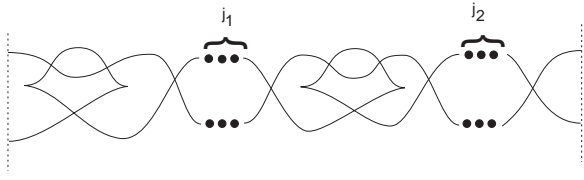}
}
\caption{The link $\mathcal L_{j_1,j_2}$.}
\label{fig:L_{j_1,j_2}}
\end{figure}

If $j_1 + j_2$ is even, then $\mathcal L_{j_1, j_2}$ is a
two-component link; else, $\mathcal L_{j_1, j_2}$ is a knot.  It is
easy to verify that $\mathcal L_{j_1, j_2}$ is topologically
isotopic to $\mathcal L_{k_1, k_2}$ if and only if $j_1 + j_2 = k_1
+ k_2$.

The next proposition shows that, for example, $\L_{1,3}$ is
distinct from $\L_{0,4}$ and from $\L_{2,2}$.

\begin{proposition}  \label{prop:even_vs_odd}
Consider $j_1, j_2$ even and $k_1, k_2$
odd with $j_1 + j_2 = k_1 + k_2$.  Then the split
Poincar\'e--Chekanov
polynomials  distinguish the links
$\mathcal L_{j_1, j_2}$ and $\mathcal L_{k_1, k_2}$.
\end{proposition}

\begin{proof} First consider the case where $k_2 > j_2$.
By ``sliding" $k_2 - j_2$ crossings around the circular base,
$\mathcal L_{k_1, k_2}$ can be rewritten as
$$(2,1, (k_1+ (k_2-j_2))^{k_2 - j_2}) \# (2,1, k_2 - (k_2 - j_2))
= (2,1,j_1^{k_2 - j_2}) \# (2,1,j_2).$$  Since $k_2 - j_2$
is odd, Theorems~\ref{thm:FlypeForm} and~\ref{thm:connect_sum} show that
$\chi^+(\lambda)[\mathcal L_{j_1, j_2}] \neq
\chi^+(\lambda)[\mathcal L_{k_1, k_2}]$.

If $k_2 < j_2$, choose $\ell$ so that
$j_2 - 2\ell < k_2 < j_2 - 2(\ell - 1)$, and
notice that
$\mathcal L_{j_1, j_2} = (2,1,(j_1 + 2\ell)^{2\ell}) \# (2,1, j_2 - 2 \ell)$.
Then by
sliding $k_2 - (j_2 - 2 \ell)$ crossings in
$\mathcal L_{k_1, k_2}$,
we have
$\L_{k_1,k_2} =  (2,1,(j_1 + 2\ell)^{k_2-(j_2 - 2\ell)}) \# (2,1, j_2
- 2\ell).$
Since $k_2 -(j_2 - 2\ell)$ is odd, the  $\chi^+$ polynomials
will again distinguish $\L_{j_1,j_2}$ and $\L_{k_1,k_2}$.
\end{proof}

\noindent
Note that, by Remark~\ref{rmk:gfconnectsum}, the generating function
polynomials
of \cite{bib:Tr} also distinguish the links in
Proposition~\ref{prop:even_vs_odd}.

%

The next proposition shows, for example, that
$\L_{0,4}$ is distinct from $\L_{2,2}$; this answers Question 1.31 in
\cite{bib:Tr}.

\begin{proposition} \label{prop:L_{k_1,k_2}_vs_L_{0,j_2}}
Consider $j_2, k_1, k_2$ with $j_2 = k_1 + k_2$ and $k_1, k_2 \neq
0$.  Then the characteristic algebras of $\mathcal L_{0, j_2}$ and
$\mathcal L_{k_1, k_2}$ are not equivalent in the sense of
\cite{bib:Ng1}, and hence the links are not Legendrian isotopic.
\end{proposition}

\begin{proof}  Let $\mathcal C$ denote the
characteristic algebra of $\mathcal L_{k_1, k_2}$,
and $\mathcal C^\prime$ that of $\mathcal L_{0,j_2}$.
It will be shown that all
elements in
$\mathcal C$ that are invertible from one side are invertible
from both sides, while in $\mathcal C^\prime$, there is an
element that is only invertible from one side.

\begin{figure}[ht]
\centerline{
\includegraphics[height=1in]{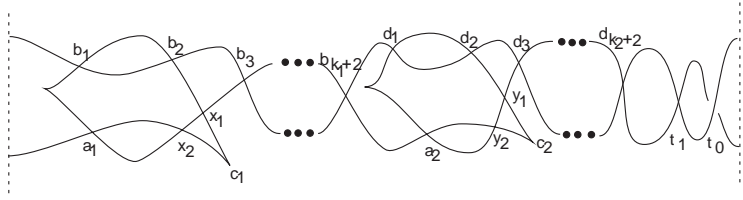}
}
\caption{A front projection of the link $\mathcal L_{k_1,k_2}$
with  its vertices labeled.}
\label{fig:L_{k_1,k_2}}
\end{figure}

First consider the Legendrian link $\mathcal L_{k_1,k_2}$.  Using
the front for $\mathcal L_{k_1,k_2}$ given in
Figure~\ref{fig:L_{k_1,k_2}}, we compute the differential on the DGA
of $\mathcal L_{k_1,k_2}$ to be
\begin{align}
\partial(c_1) &= 1 + (1 + b_2 b_1)a_1 \\
\partial(b_3) &= b_1(1 + a_1x_2)  + t_0 x_2 \label{eq:2} \\
\partial(t_1) &= t_0 \\
\partial(x_1) &= (1+b_2b_1) (1 + a_1 x_2) \label{eq:4} \\
\partial(c_2) &= 1+(1 + d_2d_1)a_2 \label{eq:5} \\
\partial(d_3) &= d_1(1 + a_2 y_2) \\
\partial(y_1) &= (1 + d_2 d_1) (1 + a_2 y_2) \label{eq:7} \\
\partial(p) &= 0, \quad  \text{ for all other vertices}.
\end{align}
In the characteristic algebra of $\mathcal L_{k_1,k_2}$, working
with the expressions given by the first four equations, we compute
that $(1+b_2 b_1)(1 + a_1 x_2) = 0$, and $b_1(1+a_1x_2)=0$,
so $1 = a_1 x_2$.  Then
$(1 + (1 + b_2 b_1)a_1) x_2 = 0$ implies
$x_2 = 1 + b_2 b_1$.  
Using this value for $x_2$, we find that
(\ref{eq:2}) transforms into
$b_1(1 + a_1(1 + b_2b_1))$, and thus (\ref{eq:4}) transforms
into
$$
(1+b_2b_1)(1+a_1(1+b_2b_1))
= 1 + a_1(1+ b_2b_1).$$
Hence, after we solve for $x_2$ and $t_0$, (\ref{eq:2})-(\ref{eq:4})
are equivalent to the relation $1 + a_1(1+ b_2b_1)$.  An analogous
argument shows
that, after solving for $y_2$, 
equations (\ref{eq:5})-(\ref{eq:7})
are equivalent to the two relations
$$1 + (1 + d_2 d_1)a_2, 1+ a_2(1 + d_2 d_1).$$
Thus   the characteristic
algebra of
$\mathcal L_{k_1, k_2}$ is
\begin{align*}
\mathcal C \simeq \mathbb{Z}/2\langle &a_1, c_1, b_1, b_2, \dots,
b_{k_1+2}, x_1, \widehat{x_2},a_2, c_2, d_1, d_2, \dots, d_{k_2+2},
y_1,
\widehat{y_2},t_1,\widehat{t_0} \rangle / \\
\langle 1 + &(1 + b_2 b_1)a_1, 1+ a_1(1 + b_2 b_1), 1 + (1 + d_2
d_1)a_2, 1+ a_2(1 + d_2 d_1) \rangle.
\end{align*}

\begin{figure}[ht]
\centerline{
\includegraphics[height=1in]{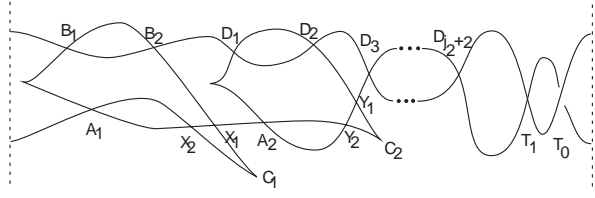}
}
\caption{A front projection of the link  $\mathcal L_{0,j_2}$
with its vertices labeled. }
\label{fig:L_{0,j_2}}
\end{figure}

Next consider the Legendrian link $\L_{0,j_2}$. Using the front in
Figure~\ref{fig:L_{0,j_2}}, we compute the differential on the DGA
of $\L_{0,j_2}$ to be
\begin{align*}
\d(C_1) &= 1 + (1+B_2B_1)A_1 + B_2T_0 \\
\d(T_1) &= T_0 \\
\d(X_1) &= (1+B_2B_1)(1+A_1X_2) + B_2T_0X_2 \\
\d(C_2) &= 1 + (1+D_2D_1)A_2 + D_2B_1(1+A_1X_2) + D_2T_0X_2 \\
\d(D_3) &= D_1(1+A_2Y_2) + B_1(1+A_1X_2)Y_2 + T_0X_2Y_2 \\
\d(Y_1) &= (1+D_2D_1)(1+A_2Y_2) + D_2B_1(1+A_1X_2)Y_2 +
D_2T_0X_2Y_2 \\
\d(P) &= 0, \quad \text{ for all other vertices.}
\end{align*}
If we quotient the resulting characteristic algebra $\C'$ by
the relations $A_2+X_2$, $B_1+1+X_2$, $B_2+1$, $D_1+1+A_1$,
$D_2+1$, $Y_2+A_1$, $C_1$, $C_2$, $X_1$, $T_0$, $T_1$, and
$D_i$ for $3\leq i\leq j_2+2$, we obtain a quotient algebra
with just two generators $A_1,X_2$ and one relation $1+A_1X_2$.
It follows that there is an element of $\C'$ (namely $A_1$) which
is invertible from the right but not from the left.
Since all elements in $\C$ which are invertible from one
side are also invertible from the other
(see also similar arguments in \cite{bib:Ng1}), we conclude
that $\C$ and $\C'$ are not equivalent.
%
\end{proof}

To use the characteristic algebra to distinguish $\L_{j_1,j_2}$ and
$\L_{k_1,k_2}$, as in
Proposition~\ref{prop:L_{k_1,k_2}_vs_L_{0,j_2}}, it is necessary
that one of $j_1,j_2,k_1,k_2$ equal $0$. If one considers the knots
$\L_{2,3}$ and $\L_{4,1}$ shown in Figure~\ref{fig:L_2_3}, then it
is not hard to verify that these knots have tamely isomorphic DGAs!
Similarly, the links $\L_{2,6}$ and $\L_{4,4}$ have isomorphic link
DGAs.  It seems likely that one could use Legendrian satellites (see
the Appendix) to distinguish these pairs, but this would  require an
involved computation with characteristic algebras. It would be
interesting to know if there was a simpler argument.


\section*{Appendix: Legendrian Satellites}
\label{sec:satellites}

In this appendix, we discuss another method to distinguish knots and
links in $\jetS$ besides generating function polynomials or the
solid-torus DGA, namely Legendrian satellites. Most of the theory
presented here either has appeared in some form previously (e.g., in
\cite{bib:Mi, bib:NgLSC}) or is part of the subject folklore, but we
include it for completeness. We then examine several examples for
which Legendrian satellites can be applied, including a pair of
solid-torus knots which cannot be distinguished by their DGAs but
can be distinguished through satellites.

Legendrian satellites are the Legendrian analogue of satellites in
the topological category. Let $L$ be an oriented Legendrian link in
$\R^3$ with one distinguished component $L_1$, and let $\tilde{L}$
be an oriented Legendrian link in $\jetS$. We give two definitions
of the Legendrian satellite $S(L,\tilde{L}) \subset \R^3$, one
abstract, one concrete.

A tubular neighborhood of $L_1$ is a solid torus contactomorphic to
$\jetS$. Thus we can embed $\tilde{L} \subset S^1 \times \R^2$ as a
Legendrian link in a tubular neighborhood of $L_1$.  Replacing the
component $L_1$ in $L$ by this new link gives $S(L,\tilde{L})$.

We can redefine $S(L,\tilde{L})$ in terms of the fronts for $L$ and
$\tilde{L}$.  First, we recall the definition of an $n$-copy from
\cite{bib:Mi}.

\begin{definition}
Given a Legendrian knot $K$ in $\R^3$, its {\it $n$-copy} is the
link consisting of $n$ copies of $K$ which differ from each other
through small perturbations in the transversal direction.  In the
front projection, the $n$-copy consists of $n$ copies of $K$,
differing from each other by small shifts in the $z$ direction.  The
$2$-copy is also known as the {\it double}.
\end{definition}

\begin{figure}
\centering{\includegraphics[width=4in]{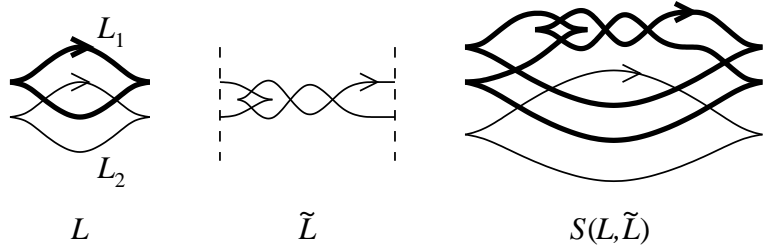}}
\caption{
Gluing a solid-torus link $\tilde{L}$ into an $\R^3$ link $L$, to
form the satellite link $S(L,\tilde{L})$.
}
\label{fig:torusgluing}
\end{figure}

Now suppose that the front of a tangle whose ends are identified to
produce $\tilde{L}$ has $2n$ endpoints, and view this front as a
Legendrian tangle in $\R^3$. Replace the front of the first
component $L_1$ of $L$ by the $n$-copy of $L_1$.  Then choose a
small segment of $L_1$ which is oriented from left to right; excise
the corresponding $n$ pieces of the $n$-copy of $L_1$, and replace
them by the front tangle for $\tilde{L}$. See
Figure~\ref{fig:torusgluing} for an illustration.

\begin{definition}
The resulting link $S(L,\tilde{L}) \subset \R^3$ is the {\it
Legendrian satellite} of $L \subset \R^3$ and $\tilde{L} \subset
S^1\times\R^2$. We give $S(L,\tilde{L})$ the orientation derived
from the orientations on $\tilde{L}$ (for the glued $n$-copy of
$L_1$) and on $L$ (for the components of $L$ besides $L_1$).
\end{definition}

\begin{proposition}
$S(L,\tilde{L})$ is a well-defined operation on Legendrian isotopy
classes; that is, if we change $L,\tilde{L}$ by Legendrian
isotopies, then $S(L,\tilde{L})$ changes by a Legendrian isotopy as
well.
\label{prop:torusgluing}
\end{proposition}

\begin{proof}
This is clear from the geometric definition of Legendrian
satellites, but we can also establish it using the concrete front
definition. We need to check that if we change one of $L,\tilde{L}$
by one of the Reidemeister moves that generate Legendrian isotopy,
then $S(L,\tilde{L})$ also changes by a Legendrian isotopy. (Note
that there is an extra set of ``Reidemeister moves'' for links in
$\jetS$, corresponding to pushing a crossing or cusp from one end of
the diagram, across the dashed lines, to the other side.) This is an
easy exercise.
\end{proof}

We now present some applications of
Proposition~\ref{prop:torusgluing} to knots and links on the solid
torus. The simplest applications glue a solid-torus knot to the
standard ``flying saucer'' unknot in $\R^3$. For instance, the top
pair of knots in Figure~\ref{fig:jetspaceknots}, glued to the
unknot, produce the Chekanov $5_2$ knots in $\R^3$; since the $5_2$
knots are not Legendrian isotopic, neither are the solid-torus
knots. Similarly, the middle pair in Figure~\ref{fig:jetspaceknots}
produce the $5_2$ knots again, while the bottom pair produce a
similar pair of $7_2$ knots. In each case,
Proposition~\ref{prop:torusgluing} implies that the solid-torus
knots are not Legendrian isotopic, recovering the results from
Section~\ref{sec:knots}.

\begin{figure}
\centering{
\includegraphics[width=4.9in]{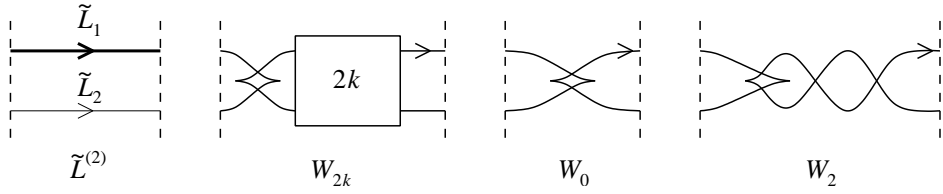}
}
\caption{
The solid-torus link $\tilde{L}^{(2)}$, and the solid-torus
Whitehead knots $W_{2k}$, $k \geq 0$, with $W_0$ and $W_2$ shown as
examples.  The box indicates $2k$ half-twists.
}
\label{fig:torus-ex}
\end{figure}

Other results can be obtained by gluing to more complicated knots or
links in $\R^3$.
  Consider the link
$\tilde{L}^{(2)}$ shown in Figure~\ref{fig:torus-ex}; we reprove the
following result from \cite{bib:Tr}, also already established in the
present paper (Corollary~\ref{cor:swap}) using the solid-torus DGA.

\begin{proposition}
Write $\tilde{L}^{(2)} = (\tilde{L}_1,\tilde{L}_2)$.  Then
$(\tilde{L}_1,\tilde{L}_2)$ is not Legendrian isotopic to
$(\tilde{L}_2,\tilde{L}_1)$.
\end{proposition}

\begin{proof}
In \cite[Prop.~4.11]{bib:Ng1}, it is proven that the double of the
usual Legendrian figure eight knot is not Legendrian isotopic to the
double with components swapped. The result now follows from
Proposition~\ref{prop:torusgluing}.
\end{proof}

Now consider the Whitehead knots $W_{2k}$ shown in
Figure~\ref{fig:torus-ex}. Each $W_{2k}$ is topologically isotopic
to its inverse (the same knot with the opposite orientation). By
contrast, we can now show the following result.

\begin{proposition}
$W_{2k}$ is not Legendrian isotopic to its inverse.
\label{prop:whiteheadknot}
\end{proposition}

\begin{proof}
Write $-W_{2k}$ for the inverse of $W_{2k}$, and let $L$ be the
double of the standard unknot in $\R^3$ (see
Figure~\ref{fig:torusgluing}). For $k=1$, it is easy to check that
$S(L,W_2)$ is precisely the oriented Whitehead link from
\cite[\S4.5]{bib:Ng1}, and that $S(L,-W_2)$ is the same link with
one component reversed. Proposition~\ref{prop:torusgluing} and
\cite[Prop.~4.12]{bib:Ng1} (which contains a typo; it should state
that $(L_6,L_7)$ and $(L_6,-L_7)$ are not Legendrian isotopic) then
imply that $W_2$ and $-W_2$ are not Legendrian isotopic.

A calculation similar to the one in the proof of
\cite[Prop.~4.12]{bib:Ng1}, omitted here, shows that $S(L,W_{2k})$
and $S(L,-W_{2k})$ are not Legendrian isotopic for arbitrary $k \geq
0$. The result follows.
\end{proof}

The solid-torus DGA even when lifted to an algebra over
$\Z[t,t^{-1}]$  fails to distinguish between $W_{2k}$ and its
inverse.  Proposition~\ref{prop:whiteheadknot} is thus a result
about solid-torus knots whose only presently known proof is via
Legendrian satellites.


\end{document}